\documentclass[10pt]{amsart}
  \usepackage{amssymb}
    \usepackage{amsfonts}
     \usepackage{textcomp}
  \usepackage[centertags]{amsmath}
    \usepackage{latexsym}
   \usepackage{amsthm}
    \usepackage{indentfirst}
 \usepackage{fancyhdr}
 \usepackage[english]{babel}
 \usepackage[english]{babel}
 \numberwithin{equation}{section}
\newtheorem{defn}{Definition}[section]
\newtheorem{lem}{Lemma}[section]
\newtheorem{thm}{Theorem}[section]
\newtheorem{rem}{Remark}[section]
\newtheorem{prop}{Proposition}[section]
\newtheorem{prf}{Proof}[section]

 \date{}
 
\begin{document}
 \title[]{\bf\Large Finsler structure for variable exponent Wasserstein space  and gradient flows}
 \author{ Aboubacar Marcos,  Ambroise Soglo} 
 \maketitle
\begin{abstract}
In this paper, we propose a variational approach based on optimal transportation to study 
the existence and unicity of solution
for a class of parabolic equations involving $q(x)$-Laplacian operator 
\begin{equation*}\label{equation variable q(x)}
 \frac{\partial \rho(t,x)}{\partial t}=div_x\left(\rho(t,x)|\nabla_x G^{'}(\rho(t,x))|^{q(x)-2}\nabla_x G^{'}(\rho(t,x))
  \right)
.\end{equation*}
The variational approach requires the setting of new tools such  as  appropiate distance on the probability space and an introduction of a Finsler metric in this space. The class of parabolic equations is derived as the flow of a gradient with respect the Finsler structure.
For $q(x)\equiv q$ constant, we recover some known results existing in the literature for the $q$-Laplacian operator. 
\end{abstract}

\section{Introduction}
The purpose of this paper is twofold:\\ Show that the Monge-Kantorovich problem associated with the Lagragian $L(x,v)=|v|^{p(x)}$ induces a distance $W_{p(.)}$ and a Finsler metric $F_ {p (.)}$ on the space of probability measures $P(\Omega)$ such that the induced distance function of $F_p (.)$ is equivalent to $W_p (.)$ where $\Omega\subset\mathbb{R}^N$ with $N\geq 1 $ is convex, and 
$\displaystyle{p(.):\Omega\longrightarrow ]1,+\infty[}$, a variable exponent function.\\
Next, establish the existence of solutions  for the following class of parabolic evolution equations involving the variable exponent
 $q(x)$ - operator.
\begin{eqnarray}\label{equation variable q(x)}
 \frac{\partial \rho(t,x)}{\partial t}&&=div_x\left(a(t,x)|\nabla_x G^{'}(\rho(t,x))|^{q(x)-2}\nabla_x G^{'}(\rho(t,x))\right), 
 \quad (t,x)\in [0,+\infty[\times\Omega\\\nonumber
 \rho(0,x)&&=\rho_0(x)\quad \mbox{in} \quad \Omega,
\end{eqnarray}
where  $\Omega\subset \mathbb{R}^N$, is a bounded domain with smooth boundary $\partial \Omega$;
$G: [0,+\infty[\rightarrow \mathbb{R}$ is a convex function of class $C^2$; 
$a: [0,+\infty[\times\Omega\rightarrow ]0,+\infty[$ is a weight function
and $q: \Omega\rightarrow ]1,+\infty[$ is a bounded measurable function
in $\Omega$, satisfying $\frac{1}{p(x)}+\frac{1}{q(x)}=1.$  Here, the initial datum $\rho_0$ is a positive
measurable function.\\
Finsler structure on the space of probability measures has been considered recently in \cite{aguey} in the case of the constant exponent (p(x) = p). Our work  generalizes the work of \cite{aguey} and moreover we derive from the Finsler metric $F_{p(.)}$,  the differential and gradient for functionals defined on the space of probability measures and then the gradient flows. Particularly, we show that the parabolic $q(x)$-Laplacian equation \eqref{equation variable q(x)} is a gradient flows of the functional 
$E(\rho)=\int_{\Omega}G(\rho)dx$ in Finslerian manifold $(P(\Omega), F_{p(.)})$ and in the Wasserstein space  $(P(\Omega),W_{p(.)}).$\\
 Equation \eqref{equation variable q(x)} presents some great interests since it is involved in the  modeling of the evolution of many nonhomogeneous materials such as   electrorheological fluids,  elastic mechanics, flow in porous media and image processing \cite{chen, Myers, Ruzicka,Zhikov}.\\ Not too many works have been devoted to the study of parabolic equations involving  variable exponent operator. In relation with our work, we can recall  the works in \cite{Huashui} and in \cite{Zhan H}.
In \cite{Zhan H}, the authors established the existence and the uniqueness  of weak solution of \eqref{equation variable q(x)} 
in the case $a(t,x)= d^{\alpha}(x)$ and $G(t)=\frac{t^2}{2}$, where $\alpha>0$ and $d(x)=dist(x,\partial \Omega)$ is the distance function from the boundary. The technique used  in their work, consist in approaching problem \eqref{equation variable q(x)}by some regularized problems under the following assumptions:  $\rho_0\in L^{\infty}(\Omega)$, 
$d(x)^{\alpha}|\nabla \rho_0|^{q^{+}}\in L^{1}(\Omega)$ and $q\equiv q(x)$ is continuous with $q^{+}:= \max_{\Omega}q(x)$.
Using similar variational approach as in \cite{Zhan H}, Huashui Zhan  \cite{Huashui} studied existence and uniqueness of solutions of \eqref{equation variable q(x)}, when $G(t)=\frac{t^2}{2}$; $a(t,x)\equiv a(x)$ and the initial datum $\rho_0$ satisfies the following assumptions
\begin{enumerate}
\item[(H1)]: $ a\in L^1(\Omega)$;  $ a^{-\frac{1}{q(x)-1}}\in L^1(\Omega)$;  $ a^{-s(x)}\in L^{1}(\Omega) $, 
with $s(x)\in ]\frac{N}{q(x)},+\infty[\cap]\frac{1}{q(x)-1},+\infty[$.\\
\item[(H2)]: $\rho_0\in L^{\infty}(\Omega)$ and $\rho_0\in W^{1,q(x)}_{a}(\Omega).$
\end{enumerate}
Here, $W^{1,q(x)}_{a}(\Omega)$ denote the  variable exponent weight Sobolev space.\\
In this paper, we propose an approach based on optimal transportation, to study existence and uniqueness of solutions of \eqref{equation variable q(x)}. From the best of our knowledge, this method is new and  requires less regularity on the initial
datum $\rho_0$ and on the variable exponent $q\equiv q(x)$ than those imposed by the authors in \cite{Huashui} and in\cite{Zhan H}  .\\
Optimal transportation method on the space of probability measures 
have be extensively used last two decades to investigate solution for parabolic partial differential equations of the form \eqref{equation variable q(x)}, when $q=q(x)$ being constant; see for instance the works in \cite{agueh,jordan}.
In \cite{jordan}, Jordan, Kinderlehrer and Otto have studied
existence of solutions of the equation \eqref{equation variable q(x)} in the particular
case where : $G(t)=\frac{t^2}{2}$, $a(t,x)\equiv 1$ and $q(x)\equiv 2$, that is the heat equation:
\begin{equation}\label{chaleur}
 \frac{\partial \rho(t,x)}{\partial t}=\Delta_x \rho(t,x),\quad \mbox{in}\quad
 [0,+\infty[\times\mathbb{R}^N.
\end{equation}
To do so, they use a descent algorithm in the probability space to construct the approximate solutions of \eqref{chaleur}. 
The descent of this algorithm is governed by the $2$- Wasserstein distance $W_2$.
In \cite{agueh},  M.Agueh used a variational approach similar as in \cite{jordan} to prove
existence of solutions for the $q$-parabolic equation
\begin{eqnarray}\label{equation variable q}
 \frac{\partial \rho(t,x)}{\partial t}&&=div_x\left(\rho(t,x)|\nabla_x G^{'}(\rho(t,x))|^{q-2}\nabla_x G^{'}(\rho(t,x))\right), 
 \quad (t,x)\in [0,+\infty[\times\Omega\\\nonumber
 \rho(0,x)&&=\rho_0(x)\quad \mbox{in} \quad \Omega
,\end{eqnarray}
 where $q$ is a constant, with $q>1$.
 Here, taking into account the works in \cite{agueh}, \cite{jordan} we extend the variational method  to the case where the  exponent is variable $ q=q(x)$  (q depends on $x$). \\
 Finally, we use the descent algorithm in the probability
space to construct the discrete solutions of \eqref{equation variable q(x)} 
and we study afterwards the convergence of our algorithm towards a weak solution of \eqref{equation variable q(x)}.
From the best of our knowledge, the solution of the Monge problem
\begin{equation}\label{monge1}
 (M):\inf_{T_{\#}\rho_1=\rho_2}\int_{\Omega}|x-T(x)|^{p(x)}\rho_1 dx
\end{equation}
is still unknown.
Howeover, if the variable exponent $p(x)$ is measurable, then the Kantorovich problem associated to \eqref{monge1}
\begin{equation}\label{kanto}
 (K):\inf_{\gamma}\left\{\int_{\mathbb{R}^N\times\mathbb{R}^N}|x-y|^{p(x)}d\gamma(x,y),\quad \gamma\in \Pi(\rho_1,\rho_2)
 \right\}
\end{equation}
admits a solution $\gamma_0$.
Furthermore, if $p(x)$ is continuous on $\Omega$, then the dual problem of \eqref{kanto}:
\begin{equation}\label{dual kanto}
 (DK) \sup_{u\oplus v \leq |x-y|^{p(x)}}\left[\int_{\Omega}u\rho_1dx+\int_{\Omega}v\rho_2 dy \right]
\end{equation}
admits a solution $(u,v)$ which satisfies the equality : $ u\oplus v=|x-y|^{p(x)}$ on the support of $\gamma_0$.\\
We then use Kantorovich's formulation to define the approximate solutions of the problem \eqref{equation variable q(x)}.\\
This paper is organised as follows.\\ Section 2 is devoted to the preliminary tools useful throughout the paper and
section 3 deals the statement of the Wasserstein distance with variable exponent. In section 4, we define the 
Finsler structure, in section 5 we give the definition of the gradient flow in Finslerian space,
and in section 6 our existence and unicity results are stated.
\section{Preliminaries}
\subsection{Main assumptions}
Throughout this work, we will assume the following:
\begin{enumerate}
 \item[$(A_1)$] $p(.):\Omega\rightarrow ]1,+\infty[$ is mesurable function such that: \\
 $1<p^{-}:=ess\inf p(x)\leq ess\sup p(x):= p^{+}<+\infty$.
 \item[$(A_2)$] $G:[0,+\infty[$ is convex function such that $G(0)=0$ and $G\in C^2(]0,+\infty[).$ 
\end{enumerate}
To prove a maximum principle for the solutions of \eqref{equation variable q(x)},\\
we assume that the initial datum $\rho_0$ is a probability density on $\Omega$ and satisfies:
\begin{enumerate}
 \item [$(A_3)$] $\rho_0 \leq M_2$, for some $0<M_2.$
 \item[$(A_4)$]
  To prove the uniqueness of solutions of \eqref{equation variable q(x)}, we assume that \\
  $\inf_{t\in[0,M_2]}G^{'}(t)>0.$
\end{enumerate}

\subsection{Lebesgue-Sobolev spaces with variable exponent}
We recall in this section some definitions and fundamental properties of the
Lebesgue and Sobolev space with variable exponents.
\begin{defn}
Let $\rho$ be a probability measure on $\Omega,$ and $p(.):\Omega\rightarrow ]1,+\infty[$ a mesurable 
function. We denote
by $L^{p(.)}_{\rho}(\Omega)$ the Lebesgue space with variable exponent defined by:
\begin{equation}
 L^{p(.)}_{\rho}(\Omega):=\left\{u:\Omega\rightarrow \mathbb{R}; \int_{\Omega}\left|\frac{u(x)}
 {\lambda}\right|^{p(x)}\rho(x)dx\leq 1, 
 \quad \mbox{for some $\lambda>0$}\right\}
,\end{equation}
with the norm
\begin{equation}\label{norme variable}
 \|u\|_{L^{p(.)}_{\rho}(\Omega)}=\inf\left\{\lambda>0, \int_{\Omega}\left|\frac{u(x)}{\lambda}
 \right|^{p(x)}\rho(x)dx\leq 1 \right\},
\end{equation} 
for all $u\in L^{p(.)}_{\rho}(\Omega).$\\
We denote by $W^{1,p(.)}_{\rho}(\Omega)$  the Sobolev space with variable exponent defined by

\begin{equation}\label{Sobolev norm}
 W^{1,p(.)}_{\rho}(\Omega):=\left\{ u\in L^{p(.)}_{\rho}(\Omega), \quad \nabla u\in
 L^{p(.)}_{\rho}(\Omega) \right\}
\end{equation}
equipped with the norm
\begin{equation}
 \|u\|_{W^{1,p(.)}_{\rho}(\Omega)}:=\|u\|_{L^{p(.)}_{\rho}(\Omega)}+\|\nabla u\|_{L^{p(.)}_{\rho}(\Omega)}
\end{equation}
\end{defn}
It is well known from the work in \cite{fan zao} that $L^{p(.)}_{\rho}(\Omega)$ and  $ W^{1,p(.)}_{\rho}(\Omega)$
are Banach spaces  respectively with the norms \eqref{norme variable} and \eqref{Sobolev norm} .\\
 We denote by $q(.):\Omega\rightarrow]1,+\infty[$ 
the conjugate function of $p(.)$ which is defined by
$$ q(x)=\frac{p(x)}{p(x)-1},\quad \mbox{for almost all $x\in\Omega$}.$$
\begin{prop}(H\"{o}lder inequality,\cite{Livre variable}).
 Let $\rho\in P(\Omega)$ and $p(.),q(.):\Omega\rightarrow]1,+\infty[$ two mesurables functions such that 
 $\frac{1}{p(x)}+\frac{1}{q(x)}=1,$ for all $x\in\Omega.$\\
 If $u\in L^{p(.)}_{\rho}(\Omega)$ and $v\in L^{q(.)}_{\rho}(\Omega)$, then:
 $$ \int_{\Omega}|u(x)v(x)|\rho(x)dx\leq (\frac{1}{p^{-}}+\frac{1}{q^{-}})\|u\|_{L^{p(.)}_{\rho}(\Omega)}
 \|v\|_{L^{q(.)}_{\rho}(\Omega)}.$$
 Furthermore, if $p_1(.),p_2(.),p_3(.):\Omega\rightarrow]1,+\infty[$ are mesurables functions such that
 $\frac{1}{p_1(x)}=\frac{1}{p_2(x)}+\frac{1}{p_3(x)},$ for almost all $x\in\Omega,$ we have
 $$ \|uv\|_{L^{p_1(.)}_{\rho}(\Omega)}\leq 2 \|u\|_{L^{p_2(.)}_{\rho}(\Omega)}
 \|v\|_{L^{p_3(.)}_{\rho}(\Omega)},$$ for 
 $u\in L^{p_2(.)}_{\rho}(\Omega)$ and $v\in L^{p_3(.)}_{\rho}(\Omega).$
 \end{prop}
\begin{prop}\label{theo1}(see \cite{fan zao}).
 Let $\rho\in P(\Omega)$ and $p_1(.), p_2(.):\Omega\longrightarrow [1,+\infty[$ two measurables functions
such that $p_1(x)\leq p_2(x)$ on $\Omega$.  Then, we have the following continious injection:
 \begin{equation} L^{p_2(.)}_{\rho}(\Omega)\hookrightarrow L^{p_1(.)}_{\rho}(\Omega).\end{equation} 
 Furthermore, $$ \|u\|_{L^{p_1(.)}_{\rho}(\Omega)}\leq 2 \|u\|_{L^{p_2(.)}_{\rho}(\Omega)}$$
\end{prop}
\begin{thm}[,\cite{fan zao}]
 Assume that the variable exponent $p(.):\Omega\rightarrow ]1;+\infty[$ satisfies $1<p^{-}\leq p(x)\leq p^{+}<+\infty$ for
 almost all $x\in \Omega$;
 with $p^{-}=ess \inf p(x)$, $p^{+}=ess \sup p(x)$.
 Then the Banach spaces $L^{p(.)}_{\rho}(\Omega)$ and $W^{1,p(.)}_{\rho}(\Omega)$ are separable, reflexive
 and uniformly convex.
\end{thm}

\subsection{BV functions}
In this section, we present some essential properties of bounded total variation functions.
\begin{defn}
 Let $\Omega \subset \mathbb{R}^N$ is be an open subset of $\mathbb{R}^N$ and
$u \in L ^ 1(\Omega)$. The total variation of u in $\Omega$ is defined by
\begin{equation}
 |D(u)|:=\sup\{\int_{\Omega}udiv(\phi)dx,\quad \phi\in C^1_c(\Omega,\mathbb{R}^N),\quad \|\phi\|_{\infty}\leq 1\}
\end{equation}
If $|D(u)|<\infty$, then  we say that u has a bounded total variation. Let denote by $BV(\Omega)$ 
the set of all functions having bounded total variation.
\end{defn}
\begin{thm}
 The BV space equipped with the standard norm
 $\|u\|_{BV}=|Du|+\|u\|_{L^1(\Omega)}$
 is a Banach space.
 Moreover, the following injection
 \begin{equation}\label{injection BV}
  BV(\Omega)\hookrightarrow L^{r}(\Omega),\quad \mbox{ is continuous } \mbox{if $1\leq r\leq \frac{N}{N-1} $}
 \end{equation}
  \qquad and compact if $  1\leq r<\frac{N}{N-1}$, $N>1$
\end{thm}

\subsection{Generalities on Finsler manifolds}
Let $M$ be a manifold, and denote by $T_{x} M$ the tangent space at $x \in M$ and by $T M := 
\bigcup_{x \in M}\{x\}\times T_{x} M$
the tangent bundle of $M$,
that is the set of all pairs $( x , v ) \in M \times T_{x} M$. A Finsler metric on $M$ is a function
$F : T M \rightarrow [ 0 , \infty)$ such that
\begin{enumerate}
 \item[ (i)] Positivity: $ F ( x , v ) > 0 $ for all $x \in M$ and $ 0\neq v \in T_{x} M$;
\item[(ii)] Positive homogeneity: $ F ( x , \lambda v ) = \lambda F ( x , v )$ for all $\lambda > 0,\quad  x \in M$ and 
$v \in T_{x} M$;
\item[(iii)] Strong convexity: $F ( x , v + w  ) \leq F ( x , v ) + F ( x , w) $ for all $x \in M$ and $v , w \in T_{ x} M$,
with equality (when $v,w\neq0$) if and
only if $v = \lambda w $ for some $\lambda > 0$.
\end{enumerate}
In the framework of a general setting condition \rm{(iii)}  in the definition of a Finsler metric on
a differentiable manifold, is replaced by the following more general condition:
\begin{enumerate}
 \item [${(iii)}^{'}$]: For all $x\in M$, and for all $w\in T_{x}M-\{0\}$, the symmetric bilinear form
 \begin{eqnarray}
  g_{x,w}: T_{x}M\times T_{x}M &&\longrightarrow \mathbb{R}\\\nonumber
           (u,v)&& \longmapsto \frac{1}{2}\frac{\partial^{2}F^{2}(x,w+tu+sv)}{\partial t\partial s}|_{t=s=0}
 \end{eqnarray}
is positive definite. 
For the convenience of our work, we use \rm{(i)}, \rm{(ii)} and \rm{(iii)} in our definition.
\end{enumerate}

\subsection{Duality mapping}
Let $X$ is be a Banach space, $X^*$ its dual and $\textless;\textgreater$ is the duality between $X^*$ and $X.$
As in \cite{mawhin}, we defined the duality mapping between $X^*$ and $X$ corresponding 
to the normalization function $\theta$ and the properties related.
\begin{defn}
 A continious function $\theta: [0,+\infty[\rightarrow [0,+\infty[$ is called a normalization function if it is strictly
 increasing, $\theta(0)=0$ and $\theta(t)\rightarrow\infty$ if $t\rightarrow\infty.$ 
\end{defn}
\begin{defn}
 Let  $\theta: [0,+\infty[\rightarrow [0,+\infty[$ be a normalization function.
 Let $X$ be is a Banach space and $X^*$ its dual.
 The duality mapping corresponding to the normalization function $\theta$ is the set values operator
 $$ J_{\theta}:X\rightarrow P_a( X^*)$$ defined by  $J_{\theta}(0)=0$ and
 
 \begin{equation}
  J_{\theta}(x):=\left\{x^*\in X^*, \quad \textless x^*, x\textgreater=\theta(\| x\|)\|x\|,\quad
  \|x^*\|=\theta(\|x\|)\right\}, \quad \mbox{ if $x\neq 0$}
 .\end{equation}
Here $P_a(X^*)$ is the set of the subsets of $X^*.$  A such duality mapping exists because of the
Hahn Banach homomorphism theorem.
\end{defn}
\begin{thm}[ Cf \cite{mawhin},  \cite{ekeland}]
Consider the normalization function $\theta: [0,+\infty[\longrightarrow [0,+\infty[$.  Assume that 
 $X$ is reflexive, locally uniformly convex. Then the following assertions hold:
 \begin{enumerate}
  \item [(i)] $Dom(J_{\theta})=X$
\item[(ii)] $J_{\theta}$ is bijective.
  \item[(iii)]Moreover  if $X^*$ is reflexive, $Card(J_{\theta}(x))=1,$ for all $x\in X.$
  
 \end{enumerate}
\end{thm}
The proof of the theorem above can be founded in \cite{mawhin}.
\section{Wasserstein metric $W_{p(.)}$}
In this section, we define a distance on the space of the probability measures
 and then study its topology. This distance induces in the probability space,
 a metric space structure in which we obtain the solutions of the equation 
 \eqref{equation variable q(x)} as a gradient flow.
\begin{defn}(Wasserstein metric $W_{p(.)}$)\\
Let
$p(.):\Omega\longrightarrow [1,+\infty[$ a measurable function and $\rho_0,\rho_1$ two probability measures on 
$\Omega.$ We define  the Wasserstein distance $W_{p(.)}$ between $\rho_0$ and $\rho_1$ by:
            \begin{equation}             
             W_{p(.)}(\rho_0,\rho_1):=\inf_{\rho(t,.)}\left\{\|v(.,.)\|_{L^{p(.)}_{\rho}([0,1]\times\Omega)}, 
                \quad \frac{\partial \rho(t,x)}{\partial t}+div_x(\rho(t,x) v(t,x))=0,\quad
                 \rho(0)=\rho_0, \quad \rho(1)=\rho_1\right\}
            ,\end{equation}
with
\begin{equation}
 \|v(.,.)\|_{L^{p(.)}_{\rho(t,.)}([0,1]\times \Omega)}: =\inf\left\{\lambda>0, 
 \quad \int_{[0,1]\times \Omega}\left|\frac{v(t,x)}{\lambda}\right|^{p(x)}\rho(t,x)dtdx \leq 1\right\}
.\end{equation}
\end{defn}
In the following, we show  that $W_{p(.)}$  is a distance on $P(\Omega)$ and we study its topology.  
\begin{lem}\label{lemmee1}
 Let $p_1(.),p_2(.):\Omega\rightarrow [1,+\infty[$ two mesurables functions such that $p_1(x)\leq p_2(x)$, 
 for almost all $x\in\Omega.$
 Then:
 $$ W_{p_1(.)}(\rho_1,\rho_2)\leq 2 W_{p_2(.)}(\rho_1,\rho_2),$$ for all $\rho_1,\rho_2\in P(\Omega).$
\end{lem}
\begin{prf}
 Let $\rho(t,.):[0,1]\rightarrow P(\Omega)$ be a arbitrary curve which is admissible for $W_{p_2(.)}(\rho_1,\rho_2)$ and
 $v_t:\mathbb{R}^N\rightarrow\mathbb{R}^N$ be a velocity field in $\mathbb{R}^N$  satisfying:
 $v_t\in L^{p_2(.)}_{\rho_t}(\Omega),$ for all $t\geq 0$ and 
 $\displaystyle{\frac{\partial \rho(t,.)}{\partial t}+div(\rho_tv_t)=0}.$ \\Then, $t\longmapsto \rho(t,.)$ is admissible for
 $W_{p_1(.)}(\rho_1,\rho_2)$. Moreover,
 $$ W_{p_1(.)}(\rho_1,\rho_2)\leq \|v_t\|_{L^{p_1(.)}_{\rho_t}([0,1]\times\Omega)}\leq 2
 \|v_t\|_{L^{p_2(.)}_{\rho_t}([0,1]\times\Omega)}$$
Then, we get
\begin{equation}
 W_{p_1(.)}(\rho_1,\rho_2)\leq 2 W_{p_2(.)}(\rho_1,\rho_2)
\end{equation}
by taking the  infimum .  $\blacksquare$
\end{prf}
By using \eqref{lemmee1}, we establish the following 
\begin{thm}
Let $p(.):\Omega\rightarrow [1,+\infty[$ be a mesurable function.
Then $W_{p(.)}$ is a distance on $P(\Omega)$.\\
Furthermore, If $\{\rho_n\} \subset P(\Omega)$ is a sequence in $P(\Omega)$ and $\rho\in P(\Omega)$;
\begin{equation}
  W_{p(.)}(\rho_n,\rho)\stackrel{n\infty}{\longrightarrow}0 \quad\Longleftrightarrow\quad 
  \rho_n\stackrel{n\infty}{\rightharpoonup}\rho \quad\mbox{narrowly}
;\end{equation}
that is, $(\rho_n)_n$ converges to $\rho$ in the metric space $(P(\Omega), W_{p(.)})$,
if and only if $(\rho_n)_n$
converges narrowly to $\rho$ in $P(\Omega).$
\end{thm}
\begin{prf}
Let $\rho_0,\rho_1\in P(\Omega)$
and $r>1$, the Monge problem
\begin{equation}
(M): \inf_{T_{\#}\rho_0=\rho_1}\left\{\int_{\Omega}\frac{|x-T(x)|^{r}}{r}\rho_0 dx\right\} 
\end{equation}
admits a unique solution $T(x)=x-|\nabla \phi|^{r-2}\nabla \phi$, where $\phi$ is $C^1$ and
$\frac{|x|^r}{r}$- convex, ie 
\begin{equation}
 \phi(x)=\inf_{y\in \mathbb{R}^N}\left\{\frac{|x-y|^r}{r}-\psi(y)\right\}
\end{equation}
 with $\psi :\mathbb{R}^N \longrightarrow \mathbb{R}$, (see \cite{aguey}).\\
Furthermore, if $t\in [0,1]$ and  $T_t=(1-t)id+tT$ is Mc Cann's interpolation, then  
the curve $\displaystyle{\rho(t,.):t\in[0,1]\longmapsto {T_t}_{\#}\rho_0}$ is the unique (constant-speed) geodesic 
joining $\rho_0$
and $\rho_1$ in the $r$-Wasserstein space $(P(\Omega),W_r)$, (see \cite{aguey}).
 and satisfies
$\displaystyle{\frac{\partial \rho}{\partial t}+div_x(\rho(T-id))=0}$ ; $\rho(0)=\rho_0$ and $\rho(1)=\rho_1.$\\
Also, we have $\displaystyle{\int_{\Omega}\left|\frac{x-T(x)}{\lambda}\right|^{p(x)}\rho_t dx<+\infty}$ 
for $\lambda= diam(\Omega).$\\
Assume that $W_{p(.)}(\rho_0,\rho_1)=0.$ We use lemma \eqref{lemmee1} and we have
\begin{equation}
 W_{1}(\rho_0,\rho_1)\leq 2 W_{p(.)}(\rho_0,\rho_1)
.\end{equation}
Then
\begin{eqnarray}\label{separer}
 W_{p(.)}(\rho_0,\rho_1)=0 & \Longrightarrow & W_{1}(\rho_0,\rho_1)=0\\\nonumber
                           & \Longrightarrow & \rho_0=\rho_1, 
,\end{eqnarray}
since $W_{1}$ is an distance on $P(\Omega)$ defined by:
\begin{equation}
W_1(\rho_0,\rho_1)=\inf_{T_{\#}\rho_0=\rho_1}\left\{\int_{\Omega}|x-T(x)|\rho_0 dx\right\}
\end{equation}
see \cite{agueh}.\\
To prove the axiom of symmetry,
let $\rho_0,\rho_1\in P(\Omega)$ and
 $\rho(t,.)$, an arbitrary curve in $P(\Omega)$ such that $\rho(0)=\rho_0$ and $\rho(1)=\rho_1$ and
 $v_t:[0,1]\times\mathbb{R}^N\rightarrow \mathbb{R}^N$ a vector field in $\mathbb{R}^N$ 
satisfying: $$v_t\in L^{p(.)}_{\rho(t,.)}([0,1]\times\Omega) \quad \mbox{and}\quad 
\frac{\partial \rho (t)}{\partial t}+ div_x(\rho(t,.)v(t))=0.$$
The reverse curve $t\mapsto\bar{\rho}(t):=\rho(1-t)$ of $t \mapsto \rho(t,.)$ satisfies:
$\bar{\rho}(0)=\rho_1$ and $\bar{\rho}(1)=\rho_0$ and the velocity fields
 $\bar{v}(t):[0,1]\times\mathbb{R}^N\rightarrow \mathbb{R}^N$ defined by $\bar{v}(t)=v(1-t)$ satisfies
$$\frac{\partial \bar{\rho} (t)}{\partial t}+ div_x(\bar{\rho}(t)\bar{v}(t))=0.$$ So, we have
$$\|\bar {v}(t)\|_{L^{p(.)}_{\bar{\rho}(t)}}([0,1]\times\Omega)=
\| {v}(t)\|_{L^{p(.)}_{{\rho}(t)}}([0,1]\times\Omega)$$ 
and next
$$W_{p(.)}(\rho_1,\rho_0)\leq \|\bar {v}(t)\|_{L^{p(.)}_{\bar{\rho}(t)}}([0,1]\times\Omega)=
\| {v}(t)\|_{L^{p(.)}_{{\rho}(t)}}([0,1]\times\Omega).$$
Since $\rho(t,.)$ is arbitrary, we get
$$W_{p(.)}(\rho_1,\rho_0)\leq W_{p(.)}(\rho_0,\rho_1)$$
 and imilarly,
$$W_{p(.)}(\rho_0,\rho_1)\leq W_{p(.)}(\rho_1,\rho_0).$$
This concludes the symmetry axiom
$$ W_{p(.)}(\rho_0,\rho_1)=W_{p(.)}(\rho_1,\rho_0).$$
Let's now show the axiom of triangular inequality.\\
Let $\rho_0$, $\rho_1$ and $\rho_2$ three elements of $P(\Omega)$ and
 $\rho_1(t), \rho_2(t)$ two arbitrary curves of $P(\Omega)$ such that:
$\rho_1(0)=\rho_0$ , $\rho_1(1)=\rho_1$ , $\rho_2(0)=\rho_1$ and $\rho_2(1)=\rho_2$. Let
$\displaystyle{ v^1_t,  v^2_t:[0,1]\times\mathbb{R}^N\longrightarrow\mathbb{R}^N}$  
two velocity fields in $\mathbb{R}^N$
satisfying respectively
$$ \int_{[0,1]\times\Omega}|\frac{v^1_t}{\alpha}|^{p(x)}\rho^1_tdtdx<+\infty,\quad \mbox{ for some $\alpha>0$},$$ 
 $$\frac{\partial \rho_1(t)}{\partial t}+div_x (\rho_1(t)v^1_t)=0$$ and 
$$\int_{[0,1]\times\Omega}|\frac{v^2_t}{\beta}|^{p(x)}\rho_2(t)dtdx<\infty,\quad \mbox{ for some $\beta>0$ },$$

$$ \frac{\partial \rho_2(t)}{\partial t}+div_x(\rho_2(t)v^2_t)=0.$$
Consider the curve $t\mapsto \bar{\rho}(t)$ defined by
\begin{eqnarray}\label{courbe translation}\left\{\begin{array}{rl}
\bar{\rho}(t)&=\rho_1(2t),\quad \mbox{if}\quad t\in[0,\frac{1}{2}]\\
\bar{\rho}(t)&=\rho_2(2t-1)\quad \mbox{if}\quad t\in[\frac{1}{2},1]
\end{array}\right.\quad
.\end{eqnarray}
We have $\bar{\rho}(0)=\rho_1(0)=\rho_0$ and $\bar{\rho}(1)=\rho_2(1)=\rho_2$ and
the velocity fields $t\mapsto \bar{v}(t)$ defined by
\begin{eqnarray}\label{courbe translation}\left\{\begin{array}{rl}
\bar{v}(t)&=v^1(2t),\quad \mbox{if}\quad t\in[0,\frac{1}{2}]\\
\bar{v}(t)&=v^2(2t-1)\quad \mbox{if}\quad t\in[\frac{1}{2},1]
\end{array}\right.\quad
,\end{eqnarray}
satisfies
$$\frac{\partial \bar{\rho}(t)}{\partial t}+div_x(\bar{\rho(t,.)}\bar{v}(t))=0.$$
Furthermore $$\|\bar{v}(t)\|_{L^{p(.)}_{\bar{\rho}(t)}}([0,1]\times\Omega)\leq\|v^1(t)\|_
{L^{p(.)}_{\rho_1(t)}}([0,1]\times\Omega)+
\|v^2(t)\|_{L^{p(.)}_{\rho_2(t)}}([0,1]\times\Omega)$$
We have that
$$ W_{p(.)}(\rho_0,\rho_2)\leq \|\bar{v}(t)\|_{L^{p(.)}_{\bar{\rho}(t)}}([0,1]\times\Omega)\leq\|v^1(t)\|_
{L^{p(.)}_{\rho_1(t)}}([0,1]\times\Omega)+
\|v^2(t)\|_{L^{p(.)}_{\rho_2(t)}}([0,1]\times\Omega);$$
since the curves $t\mapsto \rho_1(t)$ and $t\mapsto \rho_2(t)$ are arbitrary, we obtain by taking the infimum that
$$W_{p(.)}(\rho_0,\rho_2)\leq W_{p(.)}(\rho_0,\rho_1)+W_{p(.)}(\rho_1,\rho_2).$$ 
Hence, $\displaystyle{ P_{p(.)}(\Omega)=(P(\Omega),W_{p(.)})}$ is a metric space
with the Wasserstein distance $W_{p(.)}.$\\
Let's now study the topology of $ W_{p(.)}.$\\
Let $\rho\in P(\Omega)$ and $(\rho_n)_n$ be a sequence in $P(\Omega)$ such that $W_{p(.)}(\rho_n,\rho)
\stackrel{n\infty}{\longrightarrow}0$.\\
 By using \eqref{lemmee1}, we have
\begin{equation}\label{inegalitee}
 W_{1}(\rho_n,\rho)\leq 2 W_{p(.)}(\rho_n,\rho),\quad \mbox{and}\quad  W_{p(.)}(\rho_n,\rho)\leq  W_{\infty}(\rho_n,\rho)
,\end{equation}
for all $n\in \mathbb{N},$ with
\begin{equation}
 W_{\infty}(\rho_1,\rho_2)=\inf_{\gamma}\left\{ \|x-y\|_{L^{\infty}_{\gamma}(\Omega\times\Omega)},
 \quad \gamma\in \Pi(\rho_1,\rho_2)\right\}
\end{equation}
Note that, $W_{\infty}$ is Wasserstein metric defined in \cite{ambrosio}.\\
Using \eqref{inegalitee}, we have
\begin{eqnarray}
  W_{p(.)}(\rho_n,\rho)\stackrel{n\infty}{\longrightarrow}0\quad 
  & \Longrightarrow & W_{1}(\rho_n,\rho)\stackrel{n\infty}{\longrightarrow}0\\\nonumber                                                               
           & \Longrightarrow & \rho_n\stackrel{n\infty}{\rightharpoonup}\rho \quad\mbox{narrowly},\\\nonumber
\end{eqnarray}
   \mbox{because $W_{1}$ induced the narrow topology on     $P(\Omega)$.}

Inversely, if $\rho_n\stackrel{n\infty}{\rightharpoonup}\rho \quad\mbox{narrowly},$ we have, by using
\eqref{inegalitee}
\begin{eqnarray}
 W_{\infty}(\rho_n,\rho)\stackrel{n\infty}{\longrightarrow}0\quad &\Longrightarrow & 
 W_{p(.)}(\rho_n,\rho)\stackrel{n\infty}{\longrightarrow}0
.\end{eqnarray}
$\blacksquare$
\end{prf}

\section{Finsler metric on Wasserstein space ($P(\Omega),W_{p(.)}$)}
In this section, we show  that  the space of probabilities measures can be endowed with  a Finsler metric $F_p (.)$ whose associated Minkowski norm is that of the Banach quotient space $L^{p (x)}_{\rho}(\Omega)\diagup N_{\Phi}$,
with $N_{\Phi}$  the kernel of a linear and continuous application $\Phi$ defined on $L^{p (x)}_{\rho}(\Omega)$ 
with values in the space of distributions.
First of  all, we specify the tangent space of $P (\Omega)$ at the point $\rho$ as well as its topological properties.

\subsection{Tangent space $T_{\rho}P_{p(.)}(\Omega)$ and its topology}

In this part, we give a definition of the tangent space of the Wasserstein space $P_{p(.)}(\Omega)=(P(\Omega),W_{p(.)})$
 at the point $\rho $ and we study  its topological properties later.
It is known from the work of Ambrozio \cite {ambrosio} that if
 $ r> 1 $ is a real number and if $\rho(t,.):[0,1]\rightarrow P(\Omega)$ is a 
 lipschitzian curve in the $r$-Wasserstein space $(P(\Omega),W_{r})$, then there is a velocity field 
 $v(t):\mathbb{R}^N \rightarrow \mathbb{R}^N$
 in $L^{r}_{\rho(t,.)}(\Omega)$ such that: $$\frac{\partial \rho}{\partial t}+div_x(\rho(t,.)v_t)=0.$$
Using the continuous injection $L^{p+}_{\rho}(\Omega)\hookrightarrow L^{p(.)}_{\rho}(\Omega),$ 
we obtain a similar result in the  Wasserstein space $(P(\Omega),W_{p(.)})$
\begin{lem}\label{existence du flot}
 Assume that $p(.):\Omega\rightarrow ]1,\infty[$ satisfy $\rm{(A_1)}$. Let $\rho(t,.):[0,1]\rightarrow P_{p(.)}(\Omega)$ be a 
 lipschitzian curve in $(P(\Omega),W_{p^{+}})$. Then there is a velocity field  $v(t):\mathbb{R}^N \rightarrow \mathbb{R}^N$
 in $L^{p(.)}_{\rho(t,.)}(\Omega)$ such that $$\frac{\partial \rho}{\partial t}+div_x(\rho(t,.)v_t)=0.$$
\end{lem}
Using \eqref{existence du flot}, we define the tangent space as follows
\begin{defn}
 Let $p(.):\Omega\rightarrow ]1,+\infty[$  be a variable exponent and
 $\rho\in P_{p(.)}(\Omega).$ We define the tangent space $T_{\rho}P_{p(.)}(\Omega)$
 at $\rho$, of the Wasserstein space $P_{p(.)}(\Omega)=(P(\Omega), W_{p(.)})$ as :
 $$ T_{\rho}P_{p(.)}(\Omega):=\{\nu:=-div(\rho v),\quad v\in L^{p(.)}_{\rho}(\Omega)\}$$
.\end{defn}
\begin{lem}
 Assume that $p(.):\Omega\rightarrow ]1,+\infty[$ satisfy \rm{$(A_1)$}.
 Then $X=T_{\rho}P_{p(.)}(\Omega)$ endowed with 
 the norm:
 $$ \|\nu\|_{X}:=\inf\{\|v+w\|_{L^{p(.)}_{\rho}(\Omega)}:\quad w\in L^{p(.)}_{\rho}(\Omega),\quad  div(\rho w)=0 \},$$ for
 $\nu:=-div(\rho v)\in T_{\rho}P_{p(.)}(\Omega)$ is reflexive and uniformly convex Banach space.
\end{lem}
\begin{prf}
 Define \begin{eqnarray}
         \Phi: && L^{p(.)}_{\rho}(\Omega)\longrightarrow D^{'}(\Omega)\\\nonumber
               &&  v \longmapsto -div(\rho v)
        \end{eqnarray}
The map $\Phi$ is linear and continuous on $L^{p(.)}_{\rho}(\Omega)$. Then its kernel 
$$N_{\Phi}:=\{v\in L^{p(.)}_{\rho}(\Omega):\quad -div(\rho v)=0\}$$ is a closed subset of $ L^{p(.)}_{\rho}(\Omega).$
Since $ L^{p(.)}_{\rho}(\Omega)$ is be reflexive and  uniformly convex Banach space (see \cite{Livre variable}), then the
quotient space $L^{p(.)}_{\rho}(\Omega)\diagup N_{\Phi}$ is reflexive and uniformly convex Banach
space, with the norm
$$\|\bar{v}\|_{L^{p(.)}_{\rho}(\Omega)\diagup N_{\Phi}}:=\inf_{w\in N_{\Phi}}\{\|v+w\|_{L^{p(.)}_{\rho}(\Omega)}\},$$ 
for all
$\bar{v}\in L^{p(.)}_{\rho}(\Omega)\diagup N_{\Phi}.$\\
The Tangent space $ X=T_{\rho}P_{p(.)}(\Omega)$ can be identified to the image of $\Phi$  via the isomorphism:
$$L^{p(.)}_{\rho}(\Omega)\diagup N_{\Phi}\ni\bar{v}\mapsto -div(\rho v)\in T_{\rho}P_{p(.)}(\Omega).$$
Then $T_{\rho}P_{p(.)}(\Omega)$ endowed with the norm
$$ \|\nu\|_{X}:=\inf\{\|v+w\|_{L^{p(.)}_{\rho}(\Omega)}:\quad w\in L^{p(.)}_{\rho}(\Omega),\quad  div(\rho w)=0 \},$$ for
 $\nu:=-div(\rho v)\in T_{\rho}P_{p(.)}(\Omega)$ is  reflexive and uniformly convex Banach space.
 $\blacksquare$ 
 \end{prf}
As in \cite{dinca2}, we prove that the tangent space $(X=T_{\rho}P_{p(.)}(\Omega), \|.\|_{X})$ is a smooth Banach space,
when $\rho$ belong to $L^{\infty}(\Omega)$.
\begin{lem}
 Assume that $p(.):\Omega\longrightarrow ]1,+\infty[$ satisfy $\rm{(A_1)}$ and $\rho\in L^{\infty}(\Omega).$ Then the tangent
 space $(X=T_{\rho}P_{p(.)}(\Omega), \|.\|_{X})$ is smooth Banach space. And for all $\nu_1\in X-\{0\}$, 
 the G\^ateaux derivative of $\|.\|_{X}$ at $\nu_1$ is given by
 \begin{equation}
   \|.\|^{'}_{X}(\nu_1)(\nu_2):=\|.\|^{'}_{L^{p(.)}_{\rho}(\Omega)}(v_1)( v_2)
 \end{equation}
where  $v_i\in L^{p(.)}_{\rho}(\Omega)$, $i=1,2$ satisfying $\nu_i=-div(\rho v_i)$, and $v_i$
is the unique solution of the variational problem
\begin{equation}
 (P_i):\inf_{v_i\in L^{p(.)}_{\rho}(\Omega)}\left\{ \|v_i\|_{L^{p(.)}_{\rho}(\Omega)}:
 \quad \frac{\partial \rho_i(t)}{\partial t}|_{t=0}+div(\rho v_i)=0 \quad \mbox{in $\Omega$},\quad
\mbox{ $\rho v_i.\eta=0$ on $\partial\Omega$} \right\}
\end{equation}
for some curve $t\longmapsto \rho_i(t)$ in $P(\Omega)$ such that $\rho_i(0)=\rho$, and 
$\|.\|^{'}_{L^{p(.)}_{\rho}(\Omega)}(v_1)$ is the G\^ateaux derivative of the norm $ \|.\|_{L^{p(.)}_{\rho}(\Omega)} $ at
$v_1.$
\end{lem}
\begin{prf}
As in \cite{dinca2}, the Banach space $(L^{p(.)}_{\rho}(\Omega), \|.\|_{L^{p(.)}_{\rho}(\Omega)})$ is smooth if 
$\rho\in L^{\infty}(\Omega)$, and at any $u\in L^{p(.)}_{\rho}(\Omega)-\{0\}$, the G\^ateaux derivative of the norm 
$\|.\|_{L^{p(.)}_{\rho}(\Omega)}$ is given by
\begin{equation}
 \|.\|^{'}_{L^{p(.)}_{\rho}(\Omega)}(u)(h):=\frac{\displaystyle{ \int_{\Omega} 
 p(x)\frac{|u(x)|^{p(x)-1}sign(u(x))}{\|u\|^{p(x)}_{L^{p(.)}_
 {\rho}(\Omega)}} h(x)\rho dx}}{\displaystyle{ \int_{\Omega}p(x)\frac{|u(x)|^{p(x)}}{\|u\|^{p(x)+1}_{L^{p(.)}_{\rho}(\Omega)}}
 \rho dx}}.
\end{equation}
for all $h\in L^{p(.)}_{\rho}(\Omega)$.
Then, using the fact that the norm $\|.\|_{L^{p(.)}_{\rho}(\Omega)}$ is G\^ateaux differentiable, we have
\begin{eqnarray}
 \|\nu_1+t\nu_2\|&&=\inf_{w\in N_{\Phi}}\{\|v_1+tv_2+tw\|_{L^{p(.)}_{\rho}(\Omega)}\}\\\nonumber
                 &&=\|v_1\|+t\|.\|^{'}(v_1)( v_2)+\inf_{w\in N_{\Phi}}
                 \{ t\| .\|^{'}(v_1)(w)\}+0(t)
\end{eqnarray}
for $t\in \mathbb{R}$, and $0(t)$ tends to $0$ when $t$ tends to $0.$\\
Note that for $\nu_i\in X$, there exists a unique element $v_i\in L^{p(.)}_{\rho}(\Omega)$ such that $\nu_i=-div(\rho v_i)$
and
$\|.\|^{'}(v_i)(w)=0$ for all $w\in N_{\Phi}$. Therefore $v_i\in L^{p(.)}_{\rho}(\Omega)$ is
the unique solution of the variational problem:
\begin{equation}
 (P_i):\inf_{v_i\in L^{p(.)}_{\rho}(\Omega)}\left\{ \|v_i\|_{L^{p(.)}_{\rho}(\Omega)}:
 \quad \frac{\partial \rho_i(t)}{\partial t}|_{t=0}+div(\rho v_i)=0 \quad \mbox{in $\Omega$}\quad 
\mbox{and $\rho v_i.\eta=0$ on $\partial\Omega$} \right\}
\end{equation}
for some curve $t\longmapsto \rho_i(t)$ in $P(\Omega)$ such that $\rho_i(0)=\rho$.\\
Then, we have:
$$\|\nu_i\|_{X}=\| v_i\|_{L^{p(.)}_{\rho}(\Omega)}\quad \mbox{and }\quad \inf_{w\in N_{\Phi}}
                 \{ t\| .\|^{'}(v_1)(w)\}=0.$$ We conclude that:
 \begin{equation}
  \|\nu_1+t\nu_2\|_{X}=\|\nu_1\|_{X}+t\|.\|^{'}(v_1)(v_2)+0(t)
 \end{equation}
this completes our proof.  $\blacksquare$                             
\end{prf}
We denote by  $X ^* =T^*_{\rho}P_{p(.)}(\Omega)$ the dual of $X=T_{\rho}P_{p(.)}(\Omega)$, and $X^*$ is identified with
$(N_{\Phi})^{\bot}$
via the isomorphism $$(N_{\Phi})^{\bot }\ni W\mapsto L_W\in X^*$$ with
$$L_{W}(-div(\rho V)):=\int_{\Omega}WV\rho dx.$$
Here $(N_{\Phi})^{\bot }$ is the orthogonal of $N_{\Phi}$ defined by:
$$(N_{\Phi})^{\bot}:=\{W\in L^{q(x)}_{\rho}(\Omega),\quad \int_{\Omega}WV\rho dx=0,\quad \mbox{for all $V\in N_{\Phi}$}\}.$$
$(N_{\Phi})^{\bot }$ is a closed subset of $L^{q (.)}_{\rho}(\Omega)$.\\
We define a Finsler metric $F_{p(.)}$ on Wasserstein space
$(P(\Omega),W_{p(.)})$ such that, for all $\rho\in P(\Omega)$, 
 $F_{p(.)}(\rho,.)$ coincides with $X$-norm. Thus , if 
 $\displaystyle{\rho(t,.):[0,1]\longrightarrow P(\Omega)}$
is an arbitrary curve of $P_{p(.)}(\Omega)$ and if $v_t :[0,1]\times\mathbb{R}^N\longrightarrow \mathbb{R}^N$ 
is the velocity fields along of the curve $t\mapsto \rho(t,.)$, we have
\begin{equation}\label{finsler courbe}
 F_{p(.)}(\rho(t,.),-div(\rho_tv_t)):=\|-div(\rho_t v_t) \|_{X}
.\end{equation}
Using \eqref{finsler courbe}, we give the definition of the Finsler metric. 
\begin{defn}\label{finsler}[Finsler metric]
 Let $p(.):\Omega \longrightarrow ]1,+\infty[$ be a mesurable function.\\
The map
$$ F_{p(.)}:T P_{p(.)}(\Omega)\longrightarrow [0,+\infty[$$
defined on the tangent bundle $\displaystyle{ T P_{p(.)}(\Omega):=
\bigcup_{\rho\in P_{p(.)}(\Omega)}\{\rho\}\times T_{\rho}P_{p(.)}(\Omega)}$ 
of $P_{p(.)}(\Omega)$, by
 \begin{equation}
F_{p(.)}(\rho,\nu):=\|\nu\|_{X}
, \end{equation}
with $\rho\in  P_{p(.)}(\Omega)$ and $\nu\in X:= T_{\rho}P_{p(.)}(\Omega)$
is a Finsler metric on $P_{p(.)}(\Omega).$
\end{defn}
\begin{rem}
 If $\rho\in P(\Omega),$ the Minkowski norm $F_{p(.)}(\rho,.)$ can be identified
with the $T_{\rho}P_{p(.)}(\Omega)$-norm.\\
The space $P_{p(.)}(\Omega)$ endowed with the Finsler metric $F_{p(.)}$ is a Finsler manifold.
\end{rem}
\subsection{Relation between  Finsler metric  and Wasserstein distance}
Let  use the Finsler metric $F_{p(.)}$ to calculate the length of a curve of $P(\Omega).$
Indeed, if 
$\rho=\rho(t,.):[0,1]\longrightarrow P(\Omega)$ is a curve of $P(\Omega)$;
the length of $\rho(t,.)$ in Finslerian space
$(P(\Omega),F_{p(.)})$ is defined by
\begin{equation}
 L_{F_{p(.)}}(\rho(t,.)):=\int_{0}^1 F_{p(.)}(\rho(t,.),\dot{\rho}(t))dt\quad 
\quad \mbox{with $\dot{\rho}(t)=\frac{\partial \rho(t,.)}{\partial t}$}
.\end{equation}
Then, the induced distance function $d_{F_{p(.)}}$ of $F_{p(.)}$ is a distance on $P(\Omega)$ defined by
\begin{equation}
 d_{F_{p(.)}}(\rho_0,\rho_1):=\inf_{\rho(t,.)}\left\{ L_{F_{p(.)}}(\rho(t,.)),\quad  \rho(t,.):[0,1]\longrightarrow P(\Omega),
 \quad \rho(0)=\rho_0;\quad \rho(1)=\rho_1\right\}
;\end{equation}
for $\rho_0,\rho_1\in P(\Omega).$
\begin{thm}\label{prop2} Assume that $1<p^{-}\leq p(x)\leq p^{+}<+\infty.$
 For all $\rho_0,\rho_1\in P(\Omega),$ we have
 \begin{equation}
  W_{p(.)}(\rho_0,\rho_1)\leq d_{F_{p(.)}}(\rho_0,\rho_1)\lesssim W_{p(.)}^{\alpha}(\rho_0,\rho_1)
  \quad\mbox{with $\alpha=\frac{p^{-}}{p^{+}}$}
  .\end{equation}
\end{thm}
\begin{prf}
 Let $\rho_0,\rho_1\in P_{p(.)}(\Omega)$ and $\rho(t,.): [0,1]\longrightarrow P_{p(.)}(\Omega)$ is an arbitrary curve of
 $P_{p(.)}(\Omega)$ such
 that $\rho(0)=\rho_0$ and $\rho(1)=\rho_1$ and $v(t):[0,1]\times\Omega\longrightarrow \mathbb{R}^N$
a velocity field of the curve
 $t\longmapsto \rho(t,.)$ which satisfy $\frac{\partial \rho}{\partial t}+ div_x (\rho(t,.)v(t))=0$ and
 $v(t)\in L^{p(.)}_{\rho(t,.)}(\Omega)$,
 $t\geq 0.$
 We have
 \begin{equation}
  \int_{[0,1]\times\Omega}\left|\frac {v(t,x)+w}{\lambda}\right|^{p(x)}dtdx\leq 1
\quad \mbox{if $\lambda = \|v(t)+w\|_{L^{p(.)}_{\rho(t,.)}(\Omega)}$, and $-div(\rho(t,.)w)=0$}.
\end{equation}
Then
\begin{eqnarray}\label{equivalence de distance 1}
\nonumber W_{p(.)}(\rho_0,\rho_1) && \leq \|v(t,x)+w\|_{L^{p(.)}_{\rho(t,.)}([0,1]\times\Omega)}\leq
\| v(t)+w\|_{L^{p(.)}_{\rho(t,.)}(\Omega)},
\quad\mbox{for all $t\geq 0.$}\\
 \Longrightarrow W_{p(.)}(\rho_0,\rho_1) && \leq \int_{0}^1 \|-div(\rho(t,.)v_t)\|_{X} dt\leq 
 L_{F_{p(.)}}(\rho(t,.))\\\nonumber
 \Longrightarrow W_{p(.)}(\rho_0,\rho_1) && \leq d_{F_{p(.)}}(\rho_0,\rho_1)
.\end{eqnarray}
Then, if $\rho(t,.):[0,1]\longrightarrow P_{p(.)}(\Omega)$ is an curve of $P_{p(.)}(\Omega)$ such that $\rho(0)=\rho_0$
and $\rho(1)=\rho_1$
and if $v_t :\Omega\longrightarrow \mathbb{R}^N$ is a velocity field of the curve $t\longmapsto\rho(t,.)$ satisfying 
$v_t\in L^{p^{+}}_{\rho(t,.)}(\Omega)$ and $\frac{\partial \rho(t,.)}{\partial t}+ div_x(\rho(t,.)v_t)=0$, we have
by using the continuous injection
 $L^{p^{+}}_{\rho(t,.)}(\Omega)\hookrightarrow L^{p(.)}_{\rho(t,.)}(\Omega)$
\begin{eqnarray}
 d_{F_{p(.)}}(\rho_0,\rho_1) &&\leq\int_{0}^1 F_{p(.)}(\rho(t,.),\dot{\rho}(t))dt\\\nonumber
                             &&\leq\int_{0}^1\|v_t\|_{L^{p(.)}_{\rho(t,.)}(\Omega)}dt\\\nonumber
                             &&\lesssim \int_{0}^1\|v_t\|_{L^{p^{+}}_{\rho(t,.)}(\Omega)}dt
                             \\\nonumber
                             &&\lesssim \left[\int_{0}^1\|v_t\|^{p^{+}}_{L^{p^{+}}_{\rho(t,.)}(\Omega)}dt\right]
                             ^{\frac{1}{p^{+}}}             
.\end{eqnarray}
The last relationship is a consequence of the Jensen's inequality.\\
By taking the  infimum, we obtain
\begin{equation}\label{inegaliteee}
 d_{F_{p(.)}}(\rho_0,\rho_1)\lesssim W_{p^{+}}(\rho_0,\rho_1).
\end{equation}
Notice that
\begin{equation}
 W_{p^{+}}(\rho_0,\rho_1)=\inf_{T_{\#}\rho_0=\rho_1}\left[\int_{\Omega}|x-T(x)|^{p^{+}}\rho_0 dx\right]^{\frac{1}{p^{+}}}
\end{equation}
Since  $T(\Omega)\subset\Omega,$ we have
\begin{equation}\label{diametre}
 |x-T(x)|^{p^{+}}\leq (diam(\Omega))^{p^{+}-p^{-}}|x-T(x)|^{p^{-}},\mbox{ for all $x\in\Omega$}
\end{equation}
where $diam(\Omega)$ is the diameter of $\Omega$.\\
Then, using \eqref{diametre}, we get
\begin{eqnarray}
 W_{p^{+}}(\rho_0,\rho_1)&&\leq \left[\int_{\Omega}|x-T(x)|^{p^{+}}\rho_0 dx\right]^{\frac{1}{p^{+}}}\\\nonumber
                         &&\leq (diam(\Omega))^{\frac{p^{+}-p^{-}}{p^{+}}}\left[\int_{\Omega}|x-T(x)|^{p^{-}}
                         \rho_0 dx\right]^{\frac{1}{p^{+}}}
\end{eqnarray}
for all $T:\Omega\longrightarrow \Omega$, $T_{\#}\rho_0=\rho_1$ and $\int_{\Omega}|x-T(x)|^{p^{-}}\rho_0 dx<\infty.$\\
Next, by taking the infimum, we obtain
\begin{equation}\label{inegaliteeee}
 W_{p^{+}}(\rho_0,\rho_1)\leq (diam(\Omega))^{\frac{p^{+}-p^{-}}{p^{+}}} W_{p^{-}}^{\alpha}(\rho_0,\rho_1),
 \quad\mbox{with $\alpha=\frac{p^{-}}{p^{+}}$}
.\end{equation}
We use \eqref{inegalitee},\eqref{inegaliteee} and \eqref{inegaliteeee} to derive that
\begin{equation}\label{equivalence de distance 2}
 d_{F_{p(.)}}(\rho_0,\rho_1)\lesssim W_{p^{+}}(\rho_0,\rho_1)\lesssim W_{p^{-}}^{\alpha}(\rho_0,\rho_1)\lesssim 
 W_{p(.)}^{\alpha}(\rho_0,\rho_1)
.\end{equation}
Finaly, from \eqref{equivalence de distance 1} and \eqref{equivalence de distance 2}, we obtain
\begin{equation}
 W_{p(.)}(\rho_0,\rho_1)\leq d_{F_{p(.)}}(\rho_0,\rho_1)\lesssim W_{p(.)}^{\alpha}(\rho_0,\rho_1),\quad 
 \mbox{with $\alpha=\frac{p^{-}}{p^{+}}<1$}
.\end{equation}
$\blacksquare$
\end{prf}

\section{ Gradient flows in Finslerian space $(P_{p(.)}(\Omega), F_{p(.)})$}
Now, we use the Finsler metric $F_{p(.)}$ to define the gradient flows of the functional defined on the space of
probability measures.\\ 
Let $\rho\in P(\Omega)$  be a probability measure and let $p(.): \Omega\longrightarrow ]1,+\infty[$ be a mesurable 
function as in \rm{($A_1$)}. For the tangent space $X=T_{\rho}P_{p(.)}(\Omega)$ we consider the normalization 
function and the 
related duality mapping between it and its dual $X^*=T^{*}_{\rho}P(\Omega)$ as follows 

\begin{defn}
 $$ J_{\theta}:T_{\rho}P_{p(.)}(\Omega)\rightarrow P_a( X^*)$$ defined by  $J_{\theta}(0)=0$ and
 
 \begin{equation}
  J_{\theta}(\nu):=\left\{\nu^*\in X^*,  \textless \nu^*, \nu\textgreater=
  \theta(\| V\|_{L^{p(.)}_{\rho}(\Omega)})\| V\|_{L^{p(.)}_{\rho}(\Omega)},\quad
  \|\nu^*\|=\theta(\| V\|_{L^{p(.)}_{\rho}(\Omega)})=\| |V|^{p(x)-1}\|_{L^{q(.)}_{\rho}(\Omega)}\right\} 
 \end{equation}
 if $\nu\neq 0$  and with\\
  $\nu=-div(\rho V)$,  where $V\in L^{p(x)}_{\rho}(\Omega)$ is the velocity of a curve $t\longmapsto \rho(t,.)$ of
 $P(\Omega)$ at $t=0$ and the unique solution of the variational problem:
 
  \begin{equation}\label{probleme P}
 (P):\inf_{V\in L^{p(.)}_{\rho}(\Omega)}\left\{ \int_{\Omega}\frac{|V(x)|^{p(x)}}{p(x)}\rho(x)dx: 
 \quad \frac{\partial \rho(t,x)}{\partial t}|_{t=0}+div(\rho V)=0\quad 
 \mbox{in} \quad \Omega,\quad \rho V.\eta=0\quad \mbox{on $\partial \Omega$}\right\}.
\end{equation}
Here $P_a(X^*)$ is the set of the subsets of $X^*$ and $q(x)=\frac{p(x)}{p(x)-1}$
\end{defn}
Clearly  $J_{\theta}(\nu)=\left\{\nu^*\in X^*, \quad \textless \nu^*, \nu\textgreater=
  \theta(\|V\|_{L^{p(.)}_{\rho}(\Omega)} )\right \} $  \quad when $\| V\|_{L^{p(.)}_{\rho}(\Omega)} = 1$
\begin{rem}
It is worth noticing that $ J_{\theta}(\nu) $ is  well defined since such a $V$ is unique . Besides, 
the set of strictly increasing functions $\theta$ satisfying $\theta(\| V\|_{L^{p(.)}_{\rho}(\Omega)})
=\| |V|^{p(x)-1}\|_{L^{q(.)}_{\rho}(\Omega)}$ is not empty. Indeed $\theta(\| V_1\|_{L^{p(.)}_{\rho}(\Omega)}) >
\theta(\| V_2\|_{L^{p(.)}_{\rho}(\Omega)})$  when $\| V_1\|_{L^{p(.)}_{\rho}(\Omega)} >\| V_2\|_{L^{p(.)}_{\rho}(\Omega)}$
and for any $ t \in ]0,  + \infty[ $, one can find $  V'  \in  {L^{p(.)}_{\rho}(\Omega)}$ such that
$ t = \| V'\|_{L^{p(.)}_{\rho}(\Omega)}$ and then $ \theta(t) =\| |V'|^{p(x)-1}\|_{L^{q(.)}_{\rho}(\Omega)}$
\end{rem}
Since $X=T_{\rho}P_{p(.)}(\Omega)$ is reflexive, uniformly convex and smooth Banach space, the duality
map $J_{\theta}: X\rightarrow X^*$ is bijective. Then, we obtain the gradient 
of all functional defined on $P(\Omega)$  with respect to the Finsler metric $F_{p(.)}$ 
corresponding to the normalization function
$\theta: [0,+\infty[\rightarrow [0,+\infty[$
as the image by $J^{-1}_{\theta}$ 
of its differential at the point $\rho.$ 
\begin{defn}
 Let $ E : P_{p(.)}(\Omega)\longrightarrow \mathbb{R}$ is a
functional and $\rho\in P_{p(.)}(\Omega).$\\
The differential of $E$ at $\rho$ in the Finslerian space $(P_{p(.)}(\Omega), F_{p(.)})$ (if it exists) is the 
linear and bounded 
form $D_{F_{p(.)}}E(\rho)$ defined on $T_{\rho}P_{p(.)}(\Omega)$ by
\begin{equation}
 D_{F_{p(.)}}E(\rho)(\nu):=\frac{dE(\rho(t,.))}{dt}|_{t=0}
\end{equation}
where $\nu\in T_{\rho}P_{p(.)}(\Omega)$ and $\rho(t,.):[0,1]\longrightarrow P_{p(.)}(\Omega)$ is arbitrary curve
in  $P_{p(.)}(\Omega)$ satisfying
$$\rho(0)=\rho\quad\mbox{and}\quad \displaystyle{\frac{\partial\rho(t,.) }{\partial t}|_{t=0}=\nu}$$
\end{defn}
\begin{defn}
 Let $E: P_{p(.)}(\Omega)\longrightarrow \mathbb{R}$ be a functional and $\rho\in P_{p(.)}(\Omega).$\\
 The gradient of $E$ at $\rho$ 
 in the Finslerian space $(P_{p(.)}(\Omega), F_{p(.)})$ corresponding to
 the normalization function $\theta: [0,+\infty[\rightarrow [0,+\infty[$ is a unique element
 $\nabla^{\theta}_{F_{p(.)}}E(\rho)$ (if it exists) of $T_{\rho}P_{p(.)}(\Omega)$ such that
 \begin{equation}
  D_{F_{p(.)}}E(\rho)=J_{\theta}(\nabla^{\theta}_{F_{p(.)}}E(\rho))
,\end{equation}
where $J_{\theta}: T_{\rho}P_{p(.)}(\Omega)\rightarrow T^{*}_{\rho}P(\Omega)$ is the duality map corresponding to
the normalization
function $\theta:[0,+\infty[\rightarrow [0,+\infty[.$
\end{defn}
\begin{defn}
 Let $ E: P_{p(.)}(\Omega)\longrightarrow \mathbb{R}$ be a functional. 
 The curve $t\mapsto \rho(t,.)\in P(\Omega)$ is to say
 the gradient flows of $E$ in the Finslerian space $(P_{p(.)}(\Omega), F_{p(.)})$
 corresponding to the normalization function
 $\theta: [0,+\infty[\rightarrow [0,+\infty[$ if it is solution of the
 Cauchy problem
\begin{eqnarray}  
\frac{\partial \rho(t,x)}{\partial t}&& = -\nabla^{\theta}_{F_{p(.)}}E(\rho(t,x))\quad \mbox{in} 
\quad [0,+\infty[\times \Omega\\\nonumber
\rho(t=0,x)&&=\rho_0(x)\quad\mbox{in}\quad\Omega
\end{eqnarray}
\end{defn}
\begin{thm}
 Let $G:[0,+\infty[\longrightarrow \mathbb{R}$  a convex and  
 $C^2$ function. Assume that the variable exponent $p(.):\Omega \longrightarrow ]1,+\infty[$ is mesurable and satisfy
 \rm($A_1$) and $\frac{1}{p^{-}}+\frac{1}{q^{-}}\leq 1$.
 The parabolic $q(x)$-Laplacian equation
 \begin{eqnarray}\label{equation q(x)}
\frac{\partial \rho(t,x)}{\partial t}&&=div_x\{\rho(t,x)|
\nabla_x(G^{'}(\rho(t,x)))|^{q(x)-2}\nabla_x G^{'}(\rho(t,x))\}
\quad\mbox{in}\quad [0,\infty[\times \Omega\\\nonumber
\rho(t=0,x)&&=\rho_0(x)\quad\mbox{in}\quad \Omega\\\nonumber
&&\rho(t,x)|\nabla_x( G^{'}(\rho(t,x)))|^{p(x)-2}\nabla_x G^{'}(\rho(t,x)).\eta=0\quad
\mbox{on}\quad \partial \Omega
,\end{eqnarray}
 $q(x)=\frac{p(x)}{p(x)-1};$  is a gradient flows of the
 functional $E(\rho)=\int_{\Omega}G(\rho)dx$ in the Finslerian space $(P_{p(.)}(\Omega), F_{p(.)})$ corresponding to the
 normalization function $\theta,$ with $\theta(\|V\|_{L^{p(.)}_{\rho}(\Omega)})=\| |V|^{p(x)-1}\|_
 {L^{q(.)}_{\rho}(\Omega)}$ for all $V\in L^{p(.)}_{\rho}(\Omega)$.
\end{thm}
\begin{prf}
Let $\rho\in P_{p(.)}(\Omega)$, determine $\nabla^{\theta}_{F_{p(.)}}E(\rho):=-div_x(\rho V),$
 where $V\in L^{p(.)}_{\rho}(\Omega)$ is 
the unique solution of the variational
problem
 \begin{equation}
  (P):\inf_{V\in L^{p(.)}_{\rho}(\Omega)}\left\{ \int_{\Omega}\frac{|V(x)|^{p(x)}}{p(x)}\rho(x)dx: 
  \quad \frac{\partial \rho(t,x)}{\partial t}|_{t=0}=-div_x(\rho(x) V(x)),\quad
   \rho V.\eta=0 \quad \mbox{on $\partial \Omega$}\right\}
 \end{equation}
 with $\rho(t,.): [0,1]\longrightarrow P_{p(.)}(\Omega)$  a curve in $P_{p(.)}(\Omega)$ such that $\rho(0)=\rho$
 and $\frac{\partial \rho(t,.)}{\partial t}|_{t=0}=-div_x(\rho V).$\\
 We have
 \begin{equation}
  \frac{dE(\rho(t,.))}{dt}|_{t=0}=\int_{\Omega}G^{'}(\rho)\frac{\partial \rho(t,.)}{\partial t}|_{t=0}dx=
  \int_{\Omega}\rho(x)\nabla G^{'}(\rho(x)).Vdx
 \end{equation}
Using definition of $\nabla^{\theta}_{F_{p(.)}}E(\rho)$, we have
\begin{equation}
 J_{\theta}(\nabla^{\theta}_{F_{p(.)}}E(\rho))=D_{F_{p(.)}}E(\rho)
\end{equation}

Replace $V$ by  $v = \frac{V}{\theta(\| V\|_{L^{p(.)}_{\rho}(\Omega)})^ {q(x) - 1}}$, 
  then $\|v\|_{L^{p(.)}_{\rho}(\Omega)} = 1 $ and consequently 
\begin{eqnarray}\label{B}
 \textless D_{F_{p(.)}}E(\rho) , -div_x(\frac{\rho V}{\theta(\| V\|_{L^{p(.)}_{\rho}(\Omega)})^ {q(x) - 1}})\textgreater=
 \|D_{F_{p(.)}}E(\rho)\|_{X^*}  \\\nonumber
\end{eqnarray}
with $\theta(\| V\|_{L^{p(.)}_{\rho}(\Omega)}) =\| |V|^{p(x)-1}\|_{L^{q(.)}_{\rho}(\Omega)} .$\\
On the other hand, note that, if $\frac{1}{p^{-}}+\frac{1}{q^{-}}\leq 1$,
\begin{eqnarray}
 \textless  D_{F_{p(.)}}E(\rho) , -div_x(\frac{\rho V}{\theta(\| V\|_{L^{p(.)}_{\rho}(\Omega)})^ {q(x) - 1}})\textgreater &&=
 \int_{\Omega}\frac{\nabla G^{'}(\rho).V}{\theta(\| V\|_{L^{p(.)}_{\rho}(\Omega)})^ {q(x) - 1}}\rho dx \\\nonumber
&&\leq \| \nabla G^{'}(\rho)\|_{L^{q(.)}_{\rho}(\Omega)} 
\end{eqnarray}
In particular, if  $ V=|\nabla G^{'}(\rho)|^{q(x)-2}\nabla G^{'}(\rho) $, then
$\theta(\| V\|_{L^{p(.)}_{\rho}(\Omega)})=\| \nabla G^{'}(\rho)\|_{L^{q(.)}_{\rho}(\Omega)}$, \\
with $q(x)=\frac{p(x)}{p(x)-1},$  and we get
\begin{eqnarray}
 \textless D_{F_{p(.)}}E(\rho) , -div_x(\frac{\rho V}{\theta(\| V\|_{L^{p(.)}_{\rho}(\Omega)})^ {q(x) - 1}})\textgreater =
 \int_{\Omega}\frac{|\nabla G^{'}(\rho)|^{q(x)}}{\theta(\| V\|_{L^{p(.)}_{\rho}(\Omega)})^{q(x)-1}}\rho dx\\\nonumber=
 \theta(\| V\|_{L^{p(.)}_{\rho}(\Omega)})
 \int_{\Omega}|\frac{\nabla G^{'}(\rho)}{\theta(\| V\|_{L^{p(.)}_{\rho}(\Omega)})}|^{q(x)}\rho dx
\end{eqnarray}
Since $$\theta(\| V\|_{L^{p(.)}_{\rho}(\Omega)})=\| \nabla G^{'}(\rho)\|_{L^{q(.)}_{\rho}(\Omega)}, \quad
\int_{\Omega}|\frac{\nabla G^{'}(\rho)}{\theta(\| V\|_{L^{p(.)}_{\rho}(\Omega)})}|^{q(x)}\rho dx=1$$ and then
\begin{equation}\label{A}
  \textless D_{F_{p(.)}}E(\rho) , -div_x(\frac{\rho V}{\theta(\| V\|_{L^{p(.)}_{\rho}(\Omega)})^ {q(x)- 1}})
  \textgreater=\theta(\| V\|_{L^{p(.)}_{\rho}(\Omega)})=\| \nabla G^{'}(\rho)\|_{L^{q(.)}_{\rho}(\Omega)}
.\end{equation}

By using \eqref{A} and \eqref{B}, we conclude that
\begin{eqnarray}\label{C}
 \textless D_{F_{p(.)}}E(\rho) , -div_x(\frac{\rho V}{\theta(\| V\|_{L^{p(.)}_{\rho}(\Omega)})^ {q(x)- 1}})\textgreater=
 \| D_{F_{p(.)}}E(\rho)\|_{X^*}=
 \| \nabla G^{'}(\rho)\|_{L^{q(.)}_
 {\rho}(\Omega)},\\\nonumber
 \mbox{with $\theta(\| V\|_{L^{p(.)}_{\rho}(\Omega)})=\| \nabla G^{'}(\rho)\|_{L^{q(.)}_
 {\rho}(\Omega)}$ and $V=|\nabla G^{'}(\rho)|^{q(x)-2}\nabla G^{'}(\rho)$}
\end{eqnarray}
Since the map $J_{p(.)}: X\longrightarrow X^*$ is bijective, then
 $\displaystyle{\nabla_{F_{p(.)}}E(\rho)=-div_x(\rho|\nabla G^{'}(\rho)|^{q(x)-2}\nabla G^{'}(\rho))}$ is the unique
 element of $X$ which satisfy \eqref{C}.
Then the gradient flow of $E$ is the parabolic $q(x)$-Laplacian equation \eqref{equation q(x)}.

$\blacksquare$
 \end{prf}
 \begin{thm}
  The parabolic $q(x)$-Laplacian equation
\begin{eqnarray}\label{equation q(x)}
\frac{\partial \rho(t,x)}{\partial t}&&= div_x\{\rho(t,x)|\nabla_x( G^{'}(\rho(t,x)))|^{q(x)
-2}\nabla_x G^{'}(\rho(t,x))\}
\quad\mbox{in}\quad [0,\infty[\times \Omega\\\nonumber
\rho(t=0,x)&&=\rho_0(x)\quad\mbox{in}\quad \Omega\\\nonumber
&& \rho(t,x)|\nabla_x( G^{'}(\rho(t,x)))|^{q(x)-2}\nabla_x G^{'}(\rho(t,x)).\nu=0\quad\mbox{on}
\quad \partial \Omega
,\end{eqnarray}
 is the gradient flow of the functional
$E(\rho)=\int_{\Omega}G(\rho)dx$ with respect to the Wasserstein distance  $W_{p(.)}.$
 \end{thm}
\begin{prf}
 Let $t\longmapsto \rho(t,.)$  be a solution of the  parabolic equation  involving $q(x)$-Laplacian
operator \eqref{equation q(x)}. We have
 \begin{equation}
  -\frac{d E(\rho(t,.))}{dt}=\int_{\Omega}|\nabla_x G^{'}(\rho(t,x))|^{q(x)}\rho(t,x)dx
. \end{equation}
For $0<h<1$  let define\\
  $\bar{\rho}:[0,1]\longrightarrow P_{p(.)}(\Omega),$ with $\bar{\rho}(s,.)=\rho(hs+t,.)$, \\the
 curve  joining $\rho(t,.)$ and $\rho(h+t,.)$, ie such that 
$\bar{\rho}(0,.)=\rho(t,.)$ and $\bar{\rho}(1,.)=\rho(h+t,.).$ \\A velocity fields along the curve 
$s\longmapsto \bar{\rho}(s,.)$ is
$\bar{v}(s,x)=hv(hs+t,x),$ where $v(t,x)=- \|\nabla_x G^{'}(\rho(t,x))|^{q(x)-2}\nabla_x G^{'}(\rho(t,x)),$ 
the velocity fields of the curve
$t\longmapsto \rho(t,.).$ \\So, we have 
\begin{eqnarray}
\nonumber W_{p(.)}(\rho(t,.),\rho(t+h,.))&&\leq \|h v(sh+t,.)\|_{L^{p(.)}_{\rho(hs+t,.)}([0,1]\times\Omega)}
\\\nonumber && =h \|v(sh+t,.)\|_{L^{p(.)}_{\rho(hs+t,.)}([0,1]\times\Omega)}\\                       
.\end{eqnarray}
And then, we obtain
\begin{equation}\label{derive metrique}
 |\rho^{'}|(t):=\lim_{h\rightarrow 0}\frac{W_{p(.)}(\rho(t+h,.),\rho(t,.))}{h}\leq \|v(t)\|_{L^{p(.)}_{\rho(t,.)}(\Omega)}
.\end{equation}
Here, $|\rho^{'}|(t)$ is the metric derivative of the curve $t\mapsto \rho(t,.)$ at $t$ and
$v(t,x)=-|\nabla_x G^{'}(\rho(t,x))|^{q(x)-2}\nabla_x G^{'}(\rho(t,x))$ is a velocity fields along the curve
$t\longmapsto \rho(t,.).$\\
Note that if $\int_{\Omega}|v(t,x)|^{p(x)}\rho(t,x)dx=\int_{\Omega}|\nabla_x G^{'}(\rho(t,x))|^{q(x)}\rho(t,x)dx\leq 1$ then
\begin{equation}
 \|v(t)\|^{p^{+}}_{L^{p(.)}_{\rho(t,.)}(\Omega)}\leq \int_{\Omega}|\nabla_x G^{'}(\rho(t,x))|^{q(x)}\rho(t,x)dx
,\end{equation}
and if $\int_{\Omega}|v(t,x)|^{p(x)}\rho(t,x)dx=\int_{\Omega}|\nabla_x G^{'}(\rho(t,x))|^{q(x)}\rho(t,x)dx >1,$ then
\begin{equation}
 \|v(t,.)\|^{p^{-}}_{L^{p(.)}_{\rho(t,.)}(\Omega)}\leq \int_{\Omega}|\nabla_x G^{'}(\rho(t,x))|^{q(x)}\rho(t,x)dx.
\end{equation}
Write
\begin{eqnarray}\label{majoration p(x)}\left\{\begin{array}{rl}
\alpha(t)&= p^{+},\quad \mbox{if}\quad \int_{\Omega}|\nabla_x G^{'}(\rho(t,x))|^{q(x)}\rho(t,x)dx\leq 1\\
\alpha(t)&=p^{-}\quad \mbox{if}\quad \int_{\Omega}|\nabla_x G^{'}(\rho(t,x))|^{q(x)}\rho(t,x)dx >1
\end{array}\right.\quad
,\end{eqnarray}
then we have
\begin{equation}
  \|v(t,.)\|^{\alpha(t)}_{L^{p(.)}_{\rho(t,.)}(\Omega)}\leq \int_{\Omega}|\nabla_x G^{'}(\rho(t,x))|^{q(x)}\rho(t,x)dx=-
  \frac{dE(\rho(t,.))}{dt}
,\end{equation}
and recalling  \eqref{derive metrique}, we get
\begin{equation}\label{estimation 1}
 |\rho^{'}|^{\alpha(t)}(t)\leq -\frac{dE(\rho(t,.))}{dt}
.\end{equation}
Let $\rho_0,\rho_1\in P_{p(.)}(\Omega)$, $\gamma(t,.):[0,1]\longrightarrow P_{p(.)}(\Omega)$ be an arbitrary
curve such that $\gamma(0,.)=\rho_0,$ and 
$\gamma(1,.)=\rho_1$ and  $V_t(.):[0,1]\times\Omega\longrightarrow \Omega$  the velocity fields along the curve
$t\longmapsto \gamma(t)$
We have:
\begin{eqnarray}
 E(\rho_1)-E(\rho_0)&&=\int_{[0,1]}\frac{dE(\gamma(t))}{dt}dt\\\nonumber
                    &&=\int_{[0,1]\times\Omega}\nabla_x G^{'}(\gamma(t)).v(t,x)\gamma(t,x)dtdx\\\nonumber
                    &&\leq 2\|\nabla G^{'}(\gamma(t,x))\|_{L^{q(.)}_{\gamma(t,x)}([0,1]\times\Omega)}\|v(t,x)\|_
                    {L^{p(.)}_{\gamma(t,x)}([0,1]\times\Omega)}
\end{eqnarray}
Since $t\mapsto \gamma(t,.)$ is an arbitrary curve joining $\rho_0$ and $\rho_1$ in $(P(\Omega),W_{p(.)})$, we obtain:
\begin{equation}
 |E(\rho_1)-E(\rho_0)|\leq 2\|\nabla G^{'}(\gamma(t,x))\|_{L^{q(.)}_{\gamma(t,.)}([0,1]\times\Omega)}W_{p(.)}(\rho_0,\rho_1).
\end{equation}
In the particular case  we choose $\gamma$ such that $\gamma(t,.)=\rho_0$ for $t\in [0,1[$ and $\gamma(1,.)=\rho_1,$ we obtain
\begin{equation}
 \frac{|E(\rho_1)-E(\rho_0)|}{W_{p(.)}(\rho_0,\rho_1)}\leq 2\|\nabla G^{'}(\rho_0)\|_{L^{q(.)}_{\rho_0}(\Omega)}
.\end{equation}
And hence
\begin{equation}
 | \nabla_{W_{p(.)}}E(\rho_0)|:= \varlimsup_{\rho_1\rightarrow \rho_0} \frac{|E(\rho_1)-E(\rho_0)|}{W_{p(.)}(\rho_0,\rho_1)}
 \leq 2 \|\nabla G^{'}(\rho_0)\|_{L^{q(.)}_{\rho_0}(\Omega)}
.\end{equation}
Where, $|\nabla_{W_{p(.)}}E(\rho_0)|$ is the upper gradient of the functional $E(\rho)=\int_{\Omega}G(\rho)dx$ at 
$\rho_0$ with respect to the Wasserstein distance $W_{p(.)}.$\\
When we suppose 
\begin{eqnarray}\label{majoration p(x)}\left\{\begin{array}{rl}
\beta(t)&= q^{+},\quad \mbox{if}\quad \int_{\Omega}|\nabla_x G^{'}(\rho(t,x))|^{q(x)}\rho(t,x)dx\leq 1\\
\beta(t)&=q^{-}\quad \mbox{if}\quad \int_{\Omega}|\nabla_x G^{'}(\rho(t,x))|^{q(x)}\rho(t,x)dx >1
\end{array}\right.\quad
,\end{eqnarray}
We have
\begin{equation}\label{estimation 2}
 \frac{1}{2^{\beta(t)}}|\nabla_{W_{p(.)}}E(\rho(t,.))|^{\beta(t)}\leq -\frac{dE(\rho(t,.))}{dt}=\int_{\Omega}|
 \nabla G^{'}(\rho(t,x))
 |^{q(x)}\rho(t,x)dx
.\end{equation}
Recalling inequalities \eqref{estimation 1} and \eqref{estimation 2}, we conclude that  the parabolic 
$q(x)$-Laplacian equation
\eqref{equation q(x)} is a gradient flows of $E(\rho)=\int_{\Omega}G(\rho)dx$ in Wasserstein space
$(P_{p(.)}(\Omega), W_{p(.)}).$
$\blacksquare$
\end{prf}
\section{Existence and uniqueness of solution for the parabolic $q(x)$-Laplacian equation}
In this section, we prove the existence and uniqueness of solution for the class of parabolic $q(x)$-Laplacian equations
\begin{eqnarray}\label{equation q(x)}
\frac{\partial \rho(t,x)}{\partial t}&&=\displaystyle{ div_x\{\rho(t,x)|\nabla_x( G^{'}(\rho(t,x)))|^{q(x)
-2}\nabla_x G^{'}(\rho(t,x))\}}
\quad\mbox{in}\quad [0,\infty[\times \Omega\\\nonumber
\rho(t=0,x)&&=\rho_0(x)\quad\mbox{in}\quad \Omega\\\nonumber
&& \rho(t,x)|\nabla_x( G^{'}(\rho(t,x)))|^{q(x)-2}\nabla_x G^{'}(\rho(t,x)).\nu=0\quad\mbox{on}\quad \partial \Omega
,\end{eqnarray}
by using a steepest descent method in the Wasserstein space $(P(\Omega),W_{p(.)})$.
We assume that  $G:[0,+\infty[\longrightarrow \mathbb{R}$ is a convex function, which satisfies \rm{$(A_1)$},
$q(.):=\frac{p(.)}{p(.)-1}: \Omega\longrightarrow ]1,+\infty[$ is a measurable function and
$\Omega$ being an open, convex and smooth domain of $\mathbb{R}^N$, ($N\geq1$).\\
We use analogue discrete scheme as in \cite{agueh} to define a time discretization
of the problem
\eqref{equation q(x)}. 
Indeed, we fix $h>0$ to be a time step and assume that $\rho_0$ is a probability density on $\Omega.$
Define $\rho^k$, $k\in\mathbb{N}^*$ as a solution of the variational problem
\begin{equation}\label{probleme discret}
 (P^k): \inf_{\rho\in P(\Omega)}\left\{I(\rho):= \int_{\Omega}G(\rho)dx+W_{p(.)}^h(\rho^{k-1},\rho)\right\}
,\end{equation}
where
\begin{equation}
 W_{p(.)}^h(\rho^{k-1},\rho):=\inf_{\gamma\in \Pi(\rho^{k-1},\rho)}
 \left\{\int_{\Omega\times\Omega}\frac{|x-y|^{q(x)}}{h^{p(x)-1}p(x)}d\gamma(x,y)
 \right\}
.\end{equation}
Here, $\Pi(\rho^{k-1},\rho)$ is the set of all probability measures on $\Omega\times\Omega$ whose marginals
are
$\rho^{k-1}dx$ and $\rho dy.$\\
We prove in section \eqref{section3} that the sequence $(\rho^k)_k,$ satisfies the equation 
\begin{equation}\label{scheme implicite}
 \frac{\rho^k-\rho^{k-1}}{h}=div_x \left\{\rho^k|\nabla_x G^{'}(\rho^k)|^{q(x)-2}\nabla_x (G^{'}(\rho^k))\right\} + o(h)
,\end{equation}
in a weak sense, where $o(h)$ tends to $0$ when $h$ tends to $0$ and accordingly  equation \eqref{scheme implicite}
shows that the sequence $(\rho^k)_k$ is a time discretization of \eqref{equation q(x)}.\\ 
We define $\rho^h$ as follow
 \begin{eqnarray}\label{solution approcher q(x)}\left\{\begin{array}{rl}
\rho^h(t,x)&=\rho^k(x)
\quad\mbox{if}\quad t\in [hk,h(k+1)[\\
\rho^h(t=0,x)&=\rho_0(x)\quad\mbox{if}\quad t=0
\end{array}\right.\quad
\end{eqnarray}
and we show that the sequence $(\rho^h)_h$ converges weakly to $\rho(t,x)$ which solves
the parabolic $q(x)$-Laplacian equation \eqref{equation q(x)} in a weak sense.

\subsection{Euler Lagrange equation of the problem $P^k$}\label{section3}
Here, we etablish the existence and uniqueness of the solution of problem $P^k$  and  show that the sequence $(\rho^k)_k$ is a time discretization of \eqref{equation q(x)}.\\ 
\begin{prop}
Let $\rho_0$ be a probability density on $\Omega$ such that $\int_{\Omega}G(\rho_0)dx<+\infty.$
 The problem
 \begin{equation}
   (P^1): \inf_{\rho\in P(\Omega)}\left\{ I(\rho):=\int_{\Omega}G(\rho)dx+ W_{p(.)}^h(\rho_0,\rho)\right\}
   \end{equation}
 admits a unique solution $\rho^1$ and $\int_{\Omega}G(\rho^1)dx<+\infty.$
\end{prop}
\begin{prf}
 Let denote $l$ the infimum of $I$ over  $P(\Omega).$
 Show that $l$ is finite .\\
 If $\rho=\rho_0,$ then
 $\displaystyle{\int_{\Omega}G(\rho_0)dx<+\infty}$ and $W_{p(.)}^h(\rho_0,\rho_0)=0.$\\
Let $\rho$ is an probability density on $\Omega$. Since $G$ is convex, we use Jessen's inequality and obtain:
\begin{equation}
   \int_{\Omega}G(\rho)dx+ W_{p(.)}^h(\rho_0,\rho)\geq |\Omega|G(\frac{1}{|\Omega|})
.\end{equation}
We deduce that $l$ is finite.\\
Let $(\rho^n dx)_n$ be a minimizing sequence of $(P^1)$ in $P(\Omega)$.\\
Since $P(\Omega)$ is tight, then $(\rho^n dx)_n$ converges narrowly to $\rho^1 dx$ 
in $P(\Omega),$ (up to a subsequence). Since $G$ is convex and $C^1$, we have
\begin{equation}\label{convexe}
 \int_{\Omega}G(\rho^n)dx\geq \int_{\Omega}G(\rho^1)dx+\int_{\Omega}(\rho^n-\rho^1)G^{'}(\rho^1)dx
\end{equation}
By taking the limit in \eqref{convexe}, we have
\begin{equation}\label{semi continue1}
 \liminf\int_{\Omega}G(\rho^n)dx\geq \int_{\Omega}G(\rho^1)dx
\end{equation}
Let $\gamma_n$ be a solution of Kantorovich problem
\begin{equation}
 (K): \inf_{\gamma\in \Pi(\rho_0,\rho^n)}\left\{\int_{\Omega\times\Omega}\frac{|x-y|^{p(x)}}{h^{p(x)-1}p(x)}d\gamma(x,y)
 \right\}
\end{equation}
Since $P(\Omega\times\Omega)$ is tight, $(\gamma^n)_n$ converges narrowly to a probability measure
 $\gamma_1$ in $P(\Omega\times\Omega)$, (up to a subsequence), and $\gamma_1\in\Pi(\rho_0,\rho^1).$
We also have
\begin{equation}\label{semi continue2}
 \liminf \int_{\Omega\times\Omega}\frac{|x-y|^{p(x)}}{h^{p(x)-1}p(x)}d\gamma^n(x,y)
 \geq \int_{\Omega\times\Omega}\frac{|x-y|^{p(x)}}{h^{p(x)-1}p(x)}d\gamma_1(x,y)
.\end{equation}
By using \eqref{semi continue1} and \eqref{semi continue2}, we obtain
\begin{equation}
 l=\liminf\left[\int_{\Omega}G(\rho^n)dx +W_{p(.)}^h(\rho_0,\rho^n)
 \right]\geq
 \int_{\Omega}G(\rho^1)dx + W_{p(.)}^h(\rho_0,\rho^1)\geq l
.\end{equation}
Then, $\rho^1$ is a solution of $(P^1)$ and $\int_{\Omega}G(\rho^1)dx\leq \int_{\Omega}G(\rho_0)dx<\infty.$\\
We obtain uniqueness of $\rho^1$ by using the convexity of $\displaystyle{\rho\longmapsto \int_{\Omega}G(\rho)dx}$ 
and the strict convexity of the map
$\displaystyle{\rho\longmapsto W_{p(.)}^h(\rho_0,\rho)}$.
$\blacksquare$
\end{prf}
By induction, we obtain existence and uniqueness of the sequence $(\rho^k)_k$ such that $\rho^k$ is a unique solution of the
problem $(P^k).$
\begin{lem}
 Let $\rho_0$ be a probability density on $\Omega$ such that $\int_{\Omega}G(\rho_0)dx<\infty.$\\
 The Kantorovich problem
 \begin{equation}\label{Kantorovich solution}
 (K): \inf_{\gamma\in \Pi(\rho^{k-1},\rho^k)}
 \left\{\int_{\Omega\times\Omega}\frac{|x-y|^{p(x)}}{h^{p(x)-1}p(x)}d\gamma(x,y)\right\}
\end{equation}
admits a unique solution $\gamma_k$, and
$$supp\gamma_k\subset\left\{(x,y):\quad  y=x+h|\nabla_x G^{'}(\rho^k)|^{q(x)-2}\nabla_x G^{'}(\rho^k )\right\}.$$
\end{lem}
\begin{prf}
We obtain existence and uniqueness of solutions of \eqref{Kantorovich solution}  because $(x,y)\longmapsto |x-y|^{p(x)}$ 
is a Carath\'eodory function.\\
Let $\psi\in C^{\infty}_c(\Omega,\Omega)$ be a test function, and consider the flow map
$(T_{\epsilon})_{\epsilon\in\mathbb{R}}$ in $C^{\infty}_c(\mathbb{R}^n,\mathbb{R}^n)$, such that
\begin{eqnarray}\label{flot p(x)}\left\{\begin{array}{rl}
\frac{\partial T_{\epsilon}}{\partial \epsilon}&=T_{\epsilon}\circ\psi\\
T_0&= id
\end{array}\right.\quad
.\end{eqnarray}
Define: $\rho_{\epsilon}={T_{\epsilon}}_{\#}\rho^k.$\\
We have
\begin{equation}\label{deriver1}
 \frac{d}{d\epsilon}\left[\int_{\Omega}G(\rho_{\epsilon})dx\right]|_{\epsilon=0}=
 \int_{\Omega}\textless \nabla_x G^{'}(\rho^k),\psi\textgreater \rho^k dx
,\quad \mbox{see \cite{agueh}.}
\end{equation}
Let $\gamma^{\epsilon}$ be a probability measure on $\Omega\times\Omega$ defined by
\begin{equation}
 \int_{\Omega\times\Omega}\phi(x,y)d\gamma^{\epsilon}(x,y)=\int_{\Omega\times\Omega}\phi(x,T_{\epsilon}(y))
 d\gamma_k(x,y)
,\end{equation}
for all $\phi\in C^0_b(\Omega\times\Omega).$  $\gamma^{\epsilon}$
belongs to $\Pi(\rho^{k-1},\rho_{\epsilon}).$\\
Let's show that
\begin{equation}
 \frac{d}{d\epsilon}\left[\int_{\Omega\times\Omega}\frac{|x-y|^{p(x)}}{h^{p(x)-1}p(x)}
 d\gamma^{\epsilon}(x,y)\right]|_{\epsilon=0}=
 \int_{\Omega\times\Omega}\textless|\frac{x-y}{h}|^{p(x)-2}(\frac{x-y}{h}), \psi(y)\textgreater d\gamma_k(x,y)
.\end{equation}
By using the definition of $\gamma^{\epsilon},$ we have
\begin{equation}
 \int_{\Omega\times\Omega}\frac{|x-y|^{p(x)}}{h^{p(x)-1}p(x)}
 d\gamma^{\epsilon}(x,y)=\int_{\Omega\times\Omega}\frac{|x-T_{\epsilon}(y)|^{p(x)}}{h^{p(x)-1}p(x)}
 d\gamma_k(x,y)
.\end{equation}
Notice that
\begin{enumerate}
 \item[(i)] For all $\epsilon\in\mathbb{R},$ $$\int_{\Omega\times\Omega}\frac{|x-T_{\epsilon}(y)|^{p(x)}}{h^{p(x)-1}p(x)}
 d\gamma_k(x,y)<+\infty.$$ 
 \item[(ii)] For $(x,y)\in\Omega\times\Omega$ almost everywhere
 $$ \epsilon\longmapsto \frac{|x-T_{\epsilon}(y)|^{p(x)}}{h^{p(x)-1}p(x)}$$ is differentiable 
 $$\frac{d}{d\epsilon}\left[\frac{|x-T_{\epsilon}(y)|^{p(x)}}{h^{p(x)-1}p(x)} \right]=
 |\frac{x-T_{\epsilon}(y)}{h}|^{p(x)-2}(\frac{x-T_{\epsilon}(y)}{h}).\psi(T_{\epsilon}(y)).$$
Furthermore, for all $\epsilon\neq 0$
\begin{equation}
 \left|\frac{d}{d\epsilon}\left[\frac{|x-T_{\epsilon}(y)|^{p(x)}}{h^{p(x)-1}p(x)}\right]\right|\leq 
 2^{p(x)-2}|\frac{x-y}{h}|^{p(x)-1}+ \frac{2^{p(x)-2}}{h^{p(x)-1}}\|\psi\|_{\infty}
\end{equation}
 \end{enumerate}
 Recalling $(i)$, $(ii)$, and the dominated convergence theorem yield
 \begin{equation}\label{deriver2}
 \frac{d}{d\epsilon}\left[\int_{\Omega\times\Omega}\frac{|x-y|^{p(x)}}{h^{p(x)-1}p(x)}
 d\gamma^{\epsilon}(x,y)\right]|_{\epsilon=0}=
 \int_{\Omega\times\Omega}\textless|\frac{x-y}{h}|^{p(x)-2}(\frac{x-y}{h}), \psi(y)\textgreater d\gamma_k(x,y)
.\end{equation}
The solution $\rho^k$ of the problem $(P^k)$ satisfies
\begin{equation}\label{deriver3}
 \frac{d}{d\epsilon}\left[\int_{\Omega}G(\rho_{\epsilon})dx+W_{p(.)}^h(\rho^{k-1},\rho_{\epsilon}) \right]|_{\epsilon=0}=0
.\end{equation}
Note that $\gamma^{\epsilon}$ is admissible for $(P^k)$, then
\begin{equation}\label{inegaliter p(x)}
 W_{q(x)}^h(\rho^{k-1},\rho_{\epsilon})\leq \int_{\Omega}\frac{|x-y|^{p(x)}}{h^{p(x)-1}p(x)} d\gamma^{\epsilon}(x,y)
.\end{equation}
By using the previous inequality, we obtain
\begin{equation}
 I(\rho_{\epsilon})=\int_{\Omega}G(\rho_{\epsilon})dx+ W^h_{p(x)}(\rho^{k-1},\rho_{\epsilon})\leq
 \int_{\Omega}G(\rho_{\epsilon})dx+
 \int_{\Omega}\frac{|x-y|^{p(x)}}{h^{p(x)-1}p(x)} d\gamma^{\epsilon}(x,y)
.\end{equation}
So, for $\epsilon>0$, we have
\begin{eqnarray}\label{deriver4}
 \frac{I(\rho_{\epsilon})-I(\rho^k)}{\epsilon}= \displaystyle{
 \frac{\int_{\Omega}G(\rho_{\epsilon})dx-\int_{\Omega}G(\rho^k)dx}{\epsilon}}+
 \frac{W^h_{p(x)}(\rho^{k-1},\rho_{\epsilon})-W^h_{p(x)}(\rho^{k-1},\rho^k)}{\epsilon}\leq\\\nonumber
 \displaystyle{ \frac{\int_{\Omega}G(\rho_{\epsilon})dx-\int_{\Omega}G(\rho^k)dx}{\epsilon}}+
\displaystyle{\frac{\int_{\Omega\times\Omega}
 \frac{|x-T_{\epsilon}(y)|^{p(x)}}{h^{p(x)-1}p(x)}d\gamma_k(x,y)-\int_{\Omega\times\Omega}
 \frac{|x-y|^{p(x)}}{h^{p(x)-1}p(x)}d\gamma_k(x,y)}{\epsilon}}
.\end{eqnarray}
We use \eqref{deriver1}, \eqref{deriver2} , \eqref{deriver3} and \eqref{deriver4} and we tend $\epsilon$
to $0$
\begin{equation}\label{inegalite}
 0\leq \int_{\Omega}\textless \nabla_x G^{'}(\rho^k),\psi \textgreater\rho^k dx + 
 \int_{\Omega\times\Omega}\textless|\frac {x-y}{h}|^{p(x)-2}(\frac{x-y}{h}),\psi\textgreater d\gamma_k(x,y)
.\end{equation}
Changing $\psi$ by $-\psi$ in \eqref{inegalite}, we obtain the desired equality
\begin{equation}
 \int_{\Omega}\textless \nabla_x G^{'}(\rho^k),\psi \textgreater\rho^k dx +
 \int_{\Omega\times\Omega}\textless|\frac {x-y}{h}|^{p(x)-2}(\frac {x-y}{h}),\psi\textgreater d\gamma_k(x,y)=0
\end{equation}
Then
\begin{equation}\label{transport}
 y=x+h|\nabla_x G^{'}(\rho^k(x))|^{q(x)-2}\nabla_x G^{'}(\rho^k(x))\quad \gamma_k  \quad a.e
,\end{equation}
with $p(x)=\frac{q(x)}{q(x)-1}$ $\blacksquare$
\end{prf}
Now, let show that $(\rho^k)_k$ is a time discretization of \eqref{equation q(x)}.\\
Let $\psi\in C^{\infty}_c(\Omega,\mathbb{R})$ be a test function, we have
\begin{equation}\label{euler}
\int_{\Omega}(\rho^k-\rho^{k-1})\psi(x)dx=\int_{\Omega\times\Omega}(\psi(y)-\psi(x)))d\gamma_k(x,y)
\end{equation}
Using Taylor's formula
\begin{equation}\label{taylor}
 \psi(y)=\psi(x)+ (y-x).\nabla_x\psi(x)+ (y-x)^{\tau}\nabla_x^2\psi(x+\theta (y-x)).(y-x).
,\end{equation}
with  $\theta\in[0,1]$ and $(y-x)^{\tau}$ is the transpose  of $y-x$\\
We use \eqref{taylor} and \eqref{transport} in \eqref{euler}, then
\begin{eqnarray}\label{equation d'euler}
 \int_{\Omega}(\rho^k-\rho^{k-1})\psi(x)dx=-h\int_{\Omega}\textless |\nabla_x G^{'}(\rho^k(x))|^{q(x)-2}
 \nabla_x G^{'}(\rho^k(x)),\nabla_x\psi(x) \textgreater\rho^k dx-\\\nonumber
 \frac{h^2}{2}\int_{\Omega}\textless V_k^{\tau},\nabla_x^2\psi(x+\theta V_k)V_k \textgreater\rho^k dx
.\end{eqnarray}
In \eqref{equation d'euler}, $V_k:=|\nabla_x G^{'}(\rho^k(x))|^{q(x)-2}\nabla_x G^{'}(\rho^k(x)).$\\
Define $\displaystyle{ A_k(\psi)=h\int_{\Omega}\textless V_k^{\tau},\nabla_x^2\psi(x+\theta V_k)V_k \textgreater\rho^k dx}$
and show that $A_k(\psi)$ tends to $0$ when $h$ tends to $0.$\\
We have
\begin{equation}
 |A_k(\psi)|\leq h \sup_{x\in \Omega}|\nabla_x^2\psi(x)|\int_{\Omega}|\nabla_x G^{'}(\rho^k)|^{\frac{2 q(x)}{p(x)}
 }\rho^k(x)dx
\end{equation}
Since $\rho^k$ is the unique solution of $P^k$, we have
\begin{equation}
 \int_{\Omega}G(\rho^k)dx+\int_{\Omega\times\Omega}\frac{|x-y|^{p(x)}}{h^{p(x)-1}p(x)}d\gamma_k(x,y)
 \leq \int_{\Omega}G(\rho^{k-1})dx
\end{equation}
Using \eqref{transport}, we obtain
\begin{equation}\label{inegalite2}
 \int_{\Omega}G(\rho^{k-1})dx- \int_{\Omega}G(\rho^{k})dx\geq h\int_{\Omega}\frac{|\nabla_x G^{'}(\rho^k)|^{q(x)}}{p(x)}
 \rho^k dx
\end{equation}
 Using definition of $\rho^h$ and taking the sum over $k=1,..., \frac{T}{h}$ in \eqref{inegalite2}, we get
\begin{equation}\label{bornitude du gradient}
 \int_{\Omega}G(\rho_0)dx- \int_{\Omega}G(\rho^{\frac{T}{h}})dx \geq 
 \int_{[0,T]\times\Omega}\frac{|\nabla_x G^{'}(\rho^h)|^{q(x)}}{p(x)}\rho^h dtdx
\end{equation}
We use Jessen's inequality in \eqref{bornitude du gradient} and obtain:
\begin{equation}\label{bornitude du gradient 2}
 \int_{\Omega}G(\rho_0)dx- |\Omega|G(\frac{1}{|\Omega|}) \geq 
 \int_{[0, T]\times\Omega}\frac{|\nabla_x G^{'}(\rho^h)|^{q(x)}}{p(x)}\rho^h dtdx
\end{equation}

Write
\begin{equation}
 \Omega_1=\{x\in\Omega,\quad p(x)\geq 2\}\quad \mbox{and}\quad \Omega_2=\{x\in\Omega,\quad 1<p(x)< 2\}
\end{equation}
We have
\begin{equation}
 \int_{\Omega}|\nabla_x G^{'}(\rho^k)|^{\frac{2 q(x)}{p(x)}}\rho^k(x)dx=
 \int_{\Omega_1}|\nabla_x G^{'}(\rho^k)|^{\frac{2 q(x)}{p(x)}}\rho^k(x)dx+
 \int_{\Omega_2}|\nabla_x G^{'}(\rho^k)|^{\frac{2 q(x)}{p(x)}}\rho^k(x)dx
\end{equation}
Moreover,
\begin{equation}
 \int_{\Omega_1}|\nabla_x G^{'}(\rho^k)|^{\frac{2 q(x)}{p(x)}}\rho^k(x)dx\leq
  2^{q_{+}-1}\int_{\Omega_1}|\nabla_x G^{'}(\rho^k)|^{q(x)}\rho^k(x)dx + 2^{q_{+}-1}
\end{equation}
 and
 \begin{eqnarray}
 \int_{\Omega_2}|\nabla_x G^{'}(\rho^k)|^{\frac{2 q(x)}{p(x)}}\rho^k(x)dx\leq
 \int_{\Omega_2}h^{1-\frac{2}{p(x)}}|\nabla_x G^{'}(\rho^k)|^{q(x)}\rho^k(x)dx\\\nonumber
                                                                          \leq h^{1-\frac{2}{p^{-}}}
                                                                          \int_{\Omega_2}|\nabla_x G^{'}(\rho^k)|^{q(x)}
                                                                          \rho^k(x)dx
\end{eqnarray}
So, by using the previous both inequalities  and \eqref{bornitude du gradient 2}, we derive that
\begin{equation}
 |A_k(\psi)|\leq h \sup_{x\in \Omega}|\nabla_x^2\psi(x)|p^{+}\left[2^{q_{+}-1}+(2^{q_{+}-1} + h^{1-\frac{2}{p_{-}}})
 \left(\int_{\Omega}G(\rho_0)dx-
 |\Omega|G(\frac{1}{|\Omega|}\right) \right]
,\end{equation}
The inequality above  proves that $A_k(\psi)$ tends to $0$ when $h$ tends to $0$. Hence,
the sequence $(\rho^k)_k$ is a time discretization of \eqref{equation q(x)} $\blacksquare$\\
Next, let's show that the sequence $(\rho^h)_h$ converges weakly (up to a subsequence) 
to a function $\rho=\rho(t,x)$ which solves the parabolic $q(x)$-Laplacian equation \eqref{equation q(x)}.

\subsection{Convergence of the sequence $(\rho^h)_h$, weak convergence of nonlinear term \\
$div_x\{\rho^h|\nabla_x G^{'}(\rho^h)|^{p(x)-2}\nabla_x G^{'}(\rho^h)\}$}
In this section, we assume that the initial datum $\rho_0$ is a probability density which satisfies 
$\rho_0\leq M_2$ with $0<M_2$ and
we show that the function $\rho^h$ satisfies also $ \rho^h\leq M_2$, for all $h>0$
, (see \eqref{principe du maximum}). Using
\eqref{principe du maximum} and the previous results, we prove that the sequence $(H^{'}((\rho^h(t,.)))_h$ is bounded in 
$W^{1,q(x)}(\Omega),$ for all $t\geq 0$, where $H$ is a convex function such that $H^{''}(t)=tG^{''}(t)$.
Then, we use the fact that $v\longmapsto |v|^{q(x)}$ is coercive to deduce that
$(H^{'}(\rho^h(t,.)))_h$ is bounded in $W^{1,1}(\Omega)$, for all $t\geq 0.$\\
Using compactness argument on the BV spaces, we deduce that $\rho^h(t,.)$ converges strongly to 
$\rho(t,.)$ in $L^1(\Omega).$\\
Finally, we use the strong convergence of $(\rho^h(t,.))_h$ to $\rho(t,.),$ to
prove the weak convergence of the  nonlinear term
$\{div_x\{\rho^h|\nabla_x G^{'}(\rho^h)|^{q(x)-2}\nabla_x G^{'}(\rho^h)\}\}_h$ to
$div_x\{\rho|\nabla_x G^{'}(\rho)|^{q(x)-2}\nabla_x G^{'}(\rho)\}.$
\begin{prop}\label{principe du maximum}(maximum principle)\\
Assume that the initial datum $\rho_0$ satisfy $0<\rho_0\leq M_2$.\\
Then the solution $\rho^1$ of the problem $(P^1)$ satisfies: $ \rho^1\leq M_2.$
\end{prop}
\begin{prf}
 Assume by contradiction that the set $\displaystyle{E=\{y\in\Omega,\quad \rho^1(y)>M_2\}}$ 
 has a positive Lebesgue measure.\\
 Let $\gamma_1$ be a solution of Kantorovich problem
 \begin{equation}
  (K): \quad \inf_{\gamma\in \Pi(\rho_0,\rho^1)}\int_{\Omega\times\Omega}\frac{|x-y|^{p(x)}}{h^{p(x)-1}p(x)}d\gamma(x,y)
 \end{equation}
 where $\Pi(\rho_0,\rho^1)$ is the set of all probability measures on $\Omega\times\Omega$ whose marginals 
 are $\rho_0 dx$ and $\rho^1 dy.$
Denote $E^c$ the complement of $E.$
If $\gamma_1(E^c\times E)=0,$ then, we have
\begin{equation}
 |E|M<\int_{E}\rho^1(y)dy=\gamma_1(\Omega\times E)=\gamma_1(E\times E)\leq \gamma_1(E\times\Omega)=
 \int_{E}\rho_0(x)dx\leq |E|M_2
\end{equation}
which yields a contradiction. Then $\gamma_1(E^c\times E)>0.$\\
Let $\mu$ be a probability measure on $\Omega\times\Omega$ defined by
\begin{equation}
 \int_{\Omega\times\Omega}\psi(x,y)d\mu(x,y)=\int_{E^c\times E}\psi(x,y)d\gamma_1(x,y)
,\end{equation}
for all $\psi\in C_b(\Omega\times\Omega).$\\
Denote by $\mu_0$ and $\mu_1$ the marginals of $\mu$. $\mu_0$ and $\mu_1$ are absolutely continuous with 
respect to the Lebesgue measure.\\
Denote by $v_0$ and $v_1$ the respective density functions of $\mu_0$ and $\mu_1.$
We have\\
$v_0=0$ on $E$ and $v_1=0$ on $E^c.$\\
Let $\epsilon>0$ small, such that $M_2-\epsilon v_1>0$ and define: $\rho^{\epsilon}=\rho^1+\epsilon (v_0-v_1).$\\
$\rho^{\epsilon}$ is a probability density on $\Omega$ and the measure $\gamma^{\epsilon}$ defined by
\begin{equation}
 \int_{\Omega\times\Omega}\psi(x,y)d\gamma^{\epsilon}(x,y)=\int_{\Omega\times\Omega}\psi(x,y)d\gamma_1(x,y)+
 \epsilon \int_{E^c\times E}\left[\psi(x,x)-\psi(x,y)\right]d\gamma_1(x,y)
\end{equation}
belongs to $\Pi(\rho_0,\rho^{\epsilon}).$ \\
Show that $I(\rho^{\epsilon})< I(\rho^1)$.\\
The measure $\gamma^{\epsilon}$ is not necessarily a solution of Kantorovich problem
\begin{equation}
 (K): \inf_{\gamma\in \Pi(\rho_0,\rho^{\epsilon})}\int_{\Omega\times\Omega}\frac{|x-y|^{p(x)}}{h^{p(x)-1}p(x)}d\gamma(x,y).
\end{equation}
However, we have
\begin{eqnarray}\label{inegalite3}
 & I(\rho^{\epsilon})-I(\rho^1) 
\leq  \\&\left[\int_{\Omega}G(\rho^{\epsilon})dx-\int_{\Omega}G(\rho^1)dx\right] 
+
 \left[\int_{\Omega\times\Omega}\frac{|x-y|^{p(x)}}{h^{p(x)-1}p(x)}d\gamma^{\epsilon}(x,y)-
 \int_{\Omega\times\Omega}\frac{|x-y|^{p(x)}}{h^{p(x)-1}p(x)}d\gamma_1(x,y)\right]\nonumber
.\end{eqnarray}
Using definition of $\gamma^{\epsilon}$, \eqref{inegalite3} becomes
\begin{equation}
 I(\rho^{\epsilon})-I(\rho^1)\leq \left[\int_{\Omega}G(\rho^{\epsilon})dx-\int_{\Omega}G(\rho^1)dx\right]-
 \epsilon \int_{E^c\times E}\frac{|x-y|^{q(x)}}{h^{q(x)-1}q(x)}d\gamma_1(x,y)
\end{equation}
Since $G$ is $C^2$ and convex, one obtains
\begin{equation}\label{inegalite4}
 \int_{\Omega}G(\rho^{\epsilon})dx-\int_{\Omega}G(\rho^1)dx\leq -M_2\epsilon^2\int_{E^c\times E}G^{''}(\rho^1-
 \theta M_2\epsilon)d\gamma_1(x,y)<0
\end{equation}
for $\theta\in ]0,1[.$\\
Using \eqref{inegalite3} and \eqref{inegalite4}, we get $I(\rho^{\epsilon})<I(\rho^1)$, this is a contradiction because
$\rho^1$ is solution of problem $(P^1)$. We deduce that  $\rho^1\leq M_2$.
$\blacksquare$
\end{prf}
We use the definition of $\rho^h$ in \eqref{solution approcher q(x)} and then we obtain $\rho^h\leq M_2$ for all $h>0.$\\
 We use now \eqref{bornitude du gradient 2} and then
\begin{equation}
 \int_{[0,T]\times\Omega}|\nabla_x G^{'}(\rho^h)|^{q(x)}\rho^h dtdx\leq p^{+}
\left[ \int_{\Omega}G(\rho_0)dx- |\Omega| G(\frac{1}{|\Omega|})\right]
.\end{equation}
Taking into account the coercivity of  the  map:  $ v\longmapsto |v|^{q(x)} $ for any fixed $x$ when $ ess \inf q(x) > 1$, 
one can derive the following
\begin{equation}
 \int_{[0,T]\times\Omega}|\nabla_x H^{'}(\rho^h)| dx\leq K+\int_{[0,T]\times\Omega}
 |\nabla_x G^{'}(\rho^h)|^{q(x)}\rho^h dtdx\leq p^{+}
\left[ \int_{\Omega}G(\rho_0)dx- |\Omega| G(\frac{1}{|\Omega|} \right]+K.
\end{equation}
Where $K$ is a constant and $H$ is a convex function such that $H^{''}(t)=tG^{''}(t), \quad t>0.$\\
We deduce that $(\nabla_x(H^{'}(\rho^h)))_h$ is bounded on $L^{1}([0,T]\times\Omega)$.  \\
Note that $(\rho^h)_h$ is bounded in $L^{\infty}([0,T]\times\Omega),$ for $0<T<\infty$
because $\rho^h\leq M_2$ and $[0,T]\times\Omega$ is bounded when $0<T<\infty.$\\
Then $(H^{'}(\rho^h(t,.)))_h$ is bounded in $W^{1,1}(\Omega)$ for all $t\geq 0$. This implies that, up to a subsequence,
the vector valued measures $\nabla (H^{'}(\rho^h(t,.)))dtdx$ converges weakly to
a measure $\mu$ of finite  mass.\\
Hence, we have
\begin{eqnarray}
\nonumber\infty&&>\liminf_{h\rightarrow 0}\int_{[0,T]\times\Omega}|\nabla G^{'}(\rho^h)|^{q(x)}\rho^h dtdx \\
\nonumber \geq &&
\liminf_{h\rightarrow 0}
 \int_{[0,T]\times\Omega}\textless \nabla G^{'}(\rho^h), w(t,x)\textgreater\rho^hdtdx-
 \liminf_{h\rightarrow 0}\int_{[0,T]\times\Omega}|w(t,x)|^{p(x)}\rho^h dtdx\\
\nonumber &&=\int_{[0,T]\times\Omega}\textless w(t,x),\mu \textgreater dtdx-\int_{[0,T]\times\Omega}
|w(t,x)|^{p(x)}\rho(t,x)dtdx
\end{eqnarray}
for every continuous function $w.$
Consequently, $\mu$ is absolutely continuous with respect to $\rho(t,x)dtdx $ and then, there is a borel function
$K_t :[0,T]\times\Omega\rightarrow\mathbb{R}^N$, such that $\mu(dtdx)=K(t,x)\rho(t,x)dtdx.$\\
We conclude that $(H^{'}(\rho^h))$ is bounded in $BV(\Omega)$.\\
So, up to a subsequence,  there exists $\beta_t\in L^1(\Omega)$ such that $H^{'}(\rho^h(t,.))$
converge strongly to $\beta_t$ in $L^{1}(\Omega),$ \\
Since the Legendre transform  $H^*$ of $H$ is convex, then we conclude that
$\rho^h(t,.)=(H^*)^{'}(H^{'}(\rho^h))$ converge strongly to $\rho(t,.)$ in $L^1(\Omega)$.\\
 Since $G^{'}$ is continuous $G^{'}(\rho^h)$ converge strongly to $G^{'}(\rho)$ in $L^1(\Omega)$.\\
Note that,
$\displaystyle{\left||\nabla_x G^{'}(\rho^h)|^{q(x)-2}\nabla_x G^{'}(\rho^h)\right|^{p(x)}=|\nabla_x G^{'}(\rho^h)|^{q(x)}}$,
then  $(|\nabla_x G^{'}(\rho^k)|^{q(x)-2}\nabla_x G^{'}(\rho^h))_h$ is bounded in 
$L^{p(x)}([0,+\infty[\times\Omega)$.\\
 $L^{p(x)}([0,T]\times\Omega)$ being reflexive, for $0<T<\infty$,
(see \cite{fan zao}), the sequence  $(|\nabla_x G^{'}(\rho^k)|^{q(x)-2}\nabla_x G^{'}(\rho^h))_h$ converges weakly 
to some $\sigma$ in $L^{p(x)}([0,T]\times\Omega)$ up to a subsequence.\\
Arguing  as in \cite{agueh} we  show, that
$\sigma=|\nabla_x G^{'}(\rho)|^{q(x)-2}\nabla_xG^{'}(\rho)$ in  the weak sense.
\begin{thm}
 Assume that the initial datum $\rho_0$ satisfy $\rho_0\in L^{\infty}(\Omega)$ and that $G$ satisfy
 the above assumptions. If $t\longmapsto u(t)$ is a positive test function whose support is in 
 $[-T,T]$ for $0<T<\infty,$
 then
 \begin{equation}\label{non lineare}
  \lim_{h\rightarrow 0}\int_{\Omega_{T}}\textless \rho^h|\nabla_x G^{'}(\rho^h)|^{q(x)-2}\nabla_x G^{'}(\rho^h),
  \nabla_x G^{'}(\rho^h)\textgreater u(t)dtdx=\int_{\Omega_{T}}\textless \rho\sigma,\nabla_x G^{'}(\rho)\textgreater 
  u(t)dtdx
, \end{equation}
where $\Omega_{T}:=[0,T]\times\Omega$, and $\rho$ and $\sigma$ are defined above. 
Furthermore, $div_x\{\rho^h|\nabla_x G^{'}(\rho^h)|^{q(x)-2}\nabla_x G^{'}(\rho^h)\}$
converges weakly $div_x(\rho\sigma)$ in $[C^{\infty}_c(\mathbb{R}\times\Omega)]^{'}$, and
$div_x(\rho \sigma)=div_x(\rho|\nabla_x G^{'}(\rho)|^{q(x)-2}\nabla_xG^{'}(\rho) )$ in a weak sens.
\end{thm}\label{theo1}
\begin{prf}
 The proof of \eqref{non lineare} will be derived from the following three lemmas
 \begin{lem}\label{lemme1}
  For $0<T<+\infty$, we have
  \begin{equation}
   \int_{\Omega_T}\textless\rho\sigma,\nabla_x G^{'}(\rho)\textgreater u(t)dtdx  \leq 
   \liminf_{h\rightarrow 0}\int_{\Omega_T}\textless
   \rho^h|\nabla_x G^{'}(\rho^h)|^{q(x)-2}\nabla_x G^{'}(\rho^h),\nabla_x G^{'}(\rho^h)\textgreater u(t) dtdx
  \end{equation}
  with $\Omega_T:=[0,T]\times\Omega.$
 \end{lem}
\begin{prf}
Since $\rho^h$ and $u$ is positive and  $y\longmapsto |y|^{q(x)-2}y$ is monotone, we have
\begin{equation}
 \int_{\Omega_T}\rho^h\textless |\nabla_x G^{'}(\rho^h)|^{q(x)-2}\nabla_x G^{'}(\rho^h)-|\nabla_xG^{'}(\rho)|^{q(x)-2}
 \nabla_xG^{'}(\rho),\nabla_xG^{'}(\rho^h)-\nabla_xG^{'}(\rho)\textgreater u(t)dtdx\geq 0
.\end{equation}
By the previous inequality, we obtain
\begin{eqnarray}\label{lemme 1,1}
&& \int_{[0,T]\times\Omega}\rho^h\textless |\nabla_x G^{'}(\rho^h)|^{q(x)-2}\nabla_x G^{'}(\rho^h),
 \nabla_x G^{'}(\rho^h)\textgreater u(t)dtdx \geq\\\nonumber
 &&\int_{[0,T]\times\Omega}\rho^h\textless |\nabla_x G^{'}(\rho^h)|^{q(x)-2}
 \nabla_x G^{'}(\rho^h),\nabla_x G^{'}(\rho)\textgreater u(t)dtdx\\\nonumber
 &&+ \int_{[0,T]\times\Omega}\rho^h\textless |\nabla_x G^{'}(\rho)|^{q(x)-2}\nabla_x G^{'}(\rho),
 \nabla_x G^{'}(\rho^h)-\nabla_x G^{'}(\rho)\textgreater u(t)dtdx
.\end{eqnarray}
Then, using the strong convergence of $\rho^h(t,.)$ to $\rho(t,.)$, the weak convergence of \\
$(|\nabla_x G^{'}(\rho^h)|^{q(x)-2} \nabla_x G^{'}(\rho^h))_h$ to $\sigma$ and the weak convergence of
$(\nabla_x G^{'}(\rho^h))_h$ to $\nabla_x G^{'}(\rho)$, we have
\begin{equation}\label{lemme 1,2}
 \lim_{h\rightarrow0}\int_{[0,T]\times\Omega}\rho^h\textless |\nabla_x G^{'}(\rho^h)|^{q(x)-2}
 \nabla_x G^{'}(\rho^h),\nabla_x G^{'}(\rho)\textgreater u(t)dtdx=\int_{[0,T]\times\Omega}\textless \rho
 \sigma,\nabla_x G^{'}(\rho)\textgreater u(t)dtdx
.\end{equation}
Also
\begin{equation}\label{lemme 1,3}
 \lim_{h\rightarrow0}\int_{[0,T]\times\Omega}\rho^h\textless |\nabla_x G^{'}(\rho)|^{q(x)-2}\nabla_x G^{'}(\rho),
 \nabla_x G^{'}(\rho^h)-\nabla_x G^{'}(\rho)\textgreater u(t)dtdx=0
.\end{equation}
By tending $h$ to $0$ in \eqref{lemme 1,1} and  using \eqref{lemme 1,2} and \eqref{lemme 1,3},  we obtain the proof of
\eqref{lemme1}
\end{prf}
$\blacksquare$
\begin{lem}
  For $0<T<\infty$, we have
 \begin{eqnarray}\label{lemme2,0}
  &&\limsup_{h\rightarrow 0}\int_{[0,T]\times\Omega}\textless
   \rho^h|\nabla_x G^{'}(\rho^h)|^{q(x)-2}\nabla_x G^{'}(\rho^h),\nabla_x G^{'}(\rho^h)\textgreater u(t) dtdx\leq\\\nonumber
   &&
  \int_{\Omega}\left[\rho_0G^{'}(\rho_0)-G^*(G^{'}(\rho_0))\right]u(0)dx\\\nonumber
   &&+\int_{[0,T]\times\Omega}\left[\rho G^{'}(\rho)-G^*(G^{'}(\rho))\right]u^{'}(t)dtdx
 ,\end{eqnarray}
where $G^*$ is Legendre transform of $G.$
\end{lem}
\begin{prf}
Since $\rho^k$ is solution of $(P^k)$, we use \eqref{transport} and obtain
\begin{equation}
 \int_{\Omega}G(\rho^{k-1})dx-\int_{\Omega}G(\rho^k)dx\geq \frac{h}{p^{+}}
 \int_{\Omega}\textless \rho^k|\nabla_x G^{'}(\rho^k)|^{q(x)-2}\nabla_x G^{'}(\rho^k),
 \nabla_x G^{'}(\rho^k)\textgreater dx
.\end{equation}
Multiplying the previous inequality by $u\geq 0$,  we obtain after integration
\begin{equation}\label{lemme2,1}
 \sum_{k=1}^{\frac{T}{h}}\int_{t_{k-1}}^{t_k} \int_{\Omega}\left[\frac{G(\rho^{k-1})-G(\rho^k)}{h^2}\right]u(t)dtdx\geq 
 \int_{[0,T]\times\Omega}\textless \rho^h|\nabla_x G^{'}(\rho^h)|^{q(x)-2}\nabla_x G^{'}(\rho^h),
 \nabla_x G^{'}(\rho^h)\textgreater u(t) dtdx
.\end{equation}
Notice that
\begin{eqnarray}\label{lemme2,2}
 \nonumber\sum_{k=1}^{\frac{T}{h}}\int_{t_{k-1}}^{t_k} \int_{\Omega}\left[\frac{G(\rho^{k-1})-G(\rho^k)}{h^2}\right]u(t)dtdx=
 &&\int_{[0,T]\times\Omega}G(\rho^h)\left[\frac{u(t)-u(t-h^2)}{h^2}\right]dtdx\\
+ &&\frac{1}{h^2}\int_{0}^{h^2}\int_{\Omega}G(\rho^h)u(t-h^2)dtdx
.\end{eqnarray}
We tends $h$ to $0$ in \eqref{lemme2,2}, and obtain
\begin{equation}\label{lemme2,3}
 \lim_{h\rightarrow 0}\sum_{k=1}^{\frac{T}{h}}\int_{t_{k-1}}^{t_k} \int_{\Omega}\left[\frac{G(\rho^{k-1})-G(\rho^k)}{h^2}\right]u(t)dtdx=
 \int_{[0,T]\times\Omega}G(\rho)u^{'}(t)dtdx+\int_{\Omega}G(\rho_0)u(0)dx
.\end{equation}
We use \eqref{lemme2,1} and \eqref{lemme2,3}, and obtain
\begin{eqnarray}\label{lemme2,4}
 \cr&&\limsup_{h\rightarrow 0}\int_{[0,T]\times\Omega}\textless
   \rho^h|\nabla_x G^{'}(\rho^h)|^{q(x)-2}\nabla_x G^{'}(\rho^h),\nabla_x G^{'}(\rho^h)\textgreater u(t) dtdx\\\nonumber
   &&\leq 
   \left[\int_{[0,T]\times\Omega}G(\rho)u^{'}(t)dtdx+\int_{\Omega}G(\rho_0)u(0)dx\right]
\end{eqnarray}
From  the definition of $G^*$, we have $G^*(a)\geq ab-G(b)$ for all $a,b>0$ and we obtain the equality if $a=G^{'}(b).$
Then, using
$G(\rho_0)=\rho_0 G^{'}(\rho_0)-G^*(\rho_0)$ and  $G(\rho)=\rho G^{'}(\rho)- G^*(\rho)$
 in \eqref{lemme2,4}, we obtain \eqref{lemme2,0}.
 $\blacksquare$
\end{prf}
\begin{lem}
 For $0<T<\infty,$ we have
 \begin{eqnarray}\label{lemme3,0}
  \int_{[0,T]\times\Omega}\textless \rho\sigma, \nabla_x G^{'}(\rho)\textgreater u(t)dtdx\geq &&
  \int_{\Omega}\left[\rho_0G^{'}(\rho_0)-G^*(G^{'}(\rho_0))\right]u(0)dx+\\\nonumber
  &&\int_{[0,T]\times\Omega}\left[\rho G^{'}(\rho)-G^*(G^{'}(\rho))\right]u^{'}(t)dtdx
 \end{eqnarray}
\end{lem}

\begin{prf}
Define $\psi(t,x)=G^{'}(\rho(t,x))u(t),$ $\psi\in W^{1,q(x)}([0,+\infty[\times\Omega).$\\
Approximating $\psi$ by $C^{\infty}_c(\Omega)$ functions and using \eqref{equation d'euler}, we have
\begin{equation}
 \int_{\Omega}\frac{\rho^k-\rho^{k-1}}{h}\psi(t,x)dx=-\int_{\Omega}\textless 
 \rho^k|\nabla_x G^{'}(\rho^k)|^{q(x)-2}\nabla_x G^{'}(\rho^k),\nabla_x G^{'}(\rho)\textgreater u(t)dx+0(h)
,\end{equation}
where $0(h)$ tends to $0$ when $h$ tends to $0.$\\
By using the definition of  $\rho^h$, we obtain after integration
\begin{equation}\label{lemme3,6}
 \sum_{k=1}^{\frac{T}{h}}\int_{t_{k-1}}^{t_k} \int_{\Omega}\frac{\rho^k-\rho^{k-1}}{h}\psi(t,x)dtdx=
-\int_{[0,T]\times\Omega}\textless
 \rho^h |\nabla_x G^{'}(\rho^h)|^{q(x)-2}\nabla_x G^{'}(\rho^h),\nabla_x G^{'}(\rho)\textgreater u(t)dtdx+0(h) 
.\end{equation}
Also
\begin{eqnarray}\label{lemme3,1}
 \sum_{k=1}^{\frac{T}{h}}\int_{t_{k-1}}^{t_k} \int_{\Omega}\frac{\rho^k-\rho^{k-1}}{h}\psi(t,x)dtdx=
 \int_{[0,T]\times\Omega}(\rho^h-\rho)\left[\frac{\psi(t-h,x)-\psi(t,x)}{h}\right]dtdx+\\\nonumber
 -\frac{1}{h}\int_0^h\int_{\Omega}\rho^hG^{'}(\rho(t-h))u(t-h)dtdx+\\\nonumber
 \int_{[0,T]\times\Omega}\rho u(t-h)\left[\frac{G^{'}(\rho(t-h,x))-G^{'}(\rho(t,x))}{h}\right]dtdx+\\\nonumber
 \int_{[0,T]\times\Omega}\rho G^{'}(\rho)\left[\frac{u(t-h)-u(t)}{h}\right]dtdx
.\end{eqnarray}
Since $(\rho^h(t,.))_h$ converges strongly to $\rho(t,.)$
\begin{equation}\label{lemme3,2}
 \lim_{h\rightarrow0}\int_{[0,T]\times\Omega}(\rho^h-\rho)\left[\frac{\psi(t-h,x)-\psi(t,x)}{h}\right]dtdx=0
.\end{equation}
We tend $h$ to $0$ in \eqref{lemme3,1}, and using \eqref{lemme3,2}, we have
\begin{eqnarray}\label{lemme3,3}
 \nonumber \lim_{h\rightarrow0}\sum_{k=1}^{\frac{T}{h}}\int_{t_{k-1}}^{t_k} \int_{\Omega}\frac{\rho^k-\rho^{k-1}}{h}\psi(t,x)dtdx\\= 
 \lim_{h\rightarrow 0}\int_{[0,T]\times\Omega}\rho u(t-h)\left[\frac{G^{'}(\rho(t-h,x))-G^{'}(\rho(t,x))}{h}\right]dtdx+
 \\\nonumber
 -\int_{\Omega}\rho_0G^{'}(\rho_0)u(0)dx-\int_{[0,T]\times\Omega}\rho G^{'}(\rho)u^{'}(t)dtdx
.\end{eqnarray}
Since $G^*$ is convex, then
\begin{equation}
 \rho \left[ G^{'}(\rho(t-h))-G^{'}(\rho)\right]\leq  G^*(G^{'}(\rho(t-h)))-G^*(G^{'}(\rho))
.\end{equation}
Using the previous inequality, we obtain after integration
\begin{eqnarray}\label{lemme3,4}
 \int_{[0,T]\times\Omega}\rho u(t-h)\left[\frac{G^{'}(\rho(t-h))-G^{'}(\rho)}{h}\right]dtdx\leq\\
 \int_{[0,T]\times\Omega}u(t-h)\left[\frac{G^*(G^{'}(\rho(t-h)))-G^*(G^{'}(\rho))}{h}\right]dtdx\nonumber
.\end{eqnarray}
From \eqref{lemme3,4} in \eqref{lemme3,3} we obtain
\begin{eqnarray}\label{lemme3,5}
 \lim_{h\rightarrow0}\sum_{k=1}^{\frac{T}{h}}\int_{t_{k-1}}^{t_k} \int_{\Omega}\frac{\rho^k-\rho^{k-1}}{h}\psi(t,x)dtdx &&
 \leq -\int_{\Omega}\left[\rho_0G^{'}(\rho_0)-G^*(G^{'}(\rho_0))\right]u(0)dx+\\\nonumber
  &&-\int_{[0,T]\times\Omega}\left[\rho G^{'}(\rho)-G^*(G^{'}(\rho))\right]u^{'}(t)dtdx
.\end{eqnarray}
Combining \eqref{lemme3,5} and \eqref{lemme3,6} and passing to the limit, we reach\eqref{lemme3,0}.\\
 
To get the proof of \eqref{non lineare}, we use the results in the three previous lemmas .\\

Now, let show  that
$$(div_x(\rho^h|\nabla_x G^{'}(\rho^h)|^{q(x)-2}\nabla_x G^{'}(\rho^h))_h$$ converges to\\
$ div_x(\rho\sigma)=div_x(\rho |\nabla_x G^{'}(\rho)|^{q(x)-2}\nabla_x G^{'}(\rho))$ in 
$[C^{\infty}_c([0,T]\times\Omega)]^{'}$.\\

Let $\epsilon>0$ be small and $\phi\in C^{\infty}_c(\Omega)$ be a test function.
Define $\psi_{\epsilon}(t,x)=G^{'}(\rho)-\epsilon \phi(x).$\\
$\psi_{\epsilon}\in W^{1,q(x)}([0,T]\times\Omega).$ \\
We use the fact that $y\longmapsto |y|^{q(x)-2}y$ is monotone to derive
\begin{equation}
 \int_{[0,T]\times\Omega}\rho^h\textless |\nabla_x G^{'}(\rho^h)|^{q(x)-2}\nabla_xG^{'}(\rho^h)-|\nabla_x 
 \psi_{\epsilon}|^{q(x)-2}\nabla_x \psi_{\epsilon},\nabla_x G^{'}(\rho^h)-\nabla_x \psi_{\epsilon}\textgreater 
 u(t)dtdx\geq 0
.\end{equation}
Thus
\begin{eqnarray}
&&\int_{[0,T]\times\Omega}\rho^h\textless |\nabla_x G^{'}(\rho^h)|^{q(x)-2}\nabla_x G^{'}(\rho^h),\nabla_x G^{'}(\rho^h)
\textgreater u(t)dtdx \\\nonumber
&&-\int_{[0,T]\times\Omega}\rho^h\textless |\nabla_x G^{'}(\rho^h)|^{q(x)-2}\nabla_x G^{'}(\rho^h), 
\nabla_x\psi_{\epsilon}\textgreater u(t)dtdx\\\nonumber
&& -\int_{[0,T]\times\Omega}\rho^h\textless |\nabla_x \psi_{\epsilon}|^{q(x)-2}\nabla_x\psi_{\epsilon},
\nabla_x G^{'}(\rho^h)-\nabla_x\psi_{\epsilon}\textgreater u(t)dtdx\geq 0 
.\end{eqnarray}
We tends $h$ to $0$ in the previous inequality, and we use  \eqref{non lineare}, to get 
\begin{eqnarray}
 &&\int_{[0,T]\times\Omega}\rho \textless \sigma ,\nabla_x G^{'}(\rho)
\textgreater u(t)dtdx \\\nonumber
&&-\int_{[0,T]\times\Omega}\rho\textless \sigma, 
\nabla_x\psi_{\epsilon}\textgreater u(t)dtdx\\\nonumber
&& -\int_{[0,T]\times\Omega}\rho \textless |\nabla_x \psi_{\epsilon}|^{q(x)-2}\nabla_x\psi_{\epsilon},\nabla_x G^{'}(\rho)
-\nabla_x\psi_{\epsilon}\textgreater u(t)dtdx\geq 0 
.\end{eqnarray}
By using definition of $\psi_{\epsilon},$ the previous inequality becomes
\begin{equation}
 \int_{[0,T]\times\Omega}\textless \rho \sigma,\nabla_x \phi(x)\textgreater u(t)dtdx\geq 
 \int_{[0,T]\times\Omega}\rho\textless |\nabla_x \psi_{\epsilon}|^{q(x)-2}\nabla_x \psi_{\epsilon},\nabla_x\phi(x)
 \textgreater u(t) dtdx
.\end{equation}
We tends $\epsilon$ to $0$ , and we have
\begin{equation}
 \int_{[0,T]\times\Omega}\textless \rho \sigma,\nabla_x \phi(x)\textgreater u(t)dtdx\geq 
 \int_{[0,T]\times\Omega}\rho\textless |\nabla_x G^{'}(\rho)|^{q(x)-2}\nabla_x G^{'}(\rho),\nabla_x\phi(x)
 \textgreater u(t) dtdx
.\end{equation}
Replacing $\phi$ by $-\phi$ in the  previous inequality, we obtain the equality:
\begin{equation}
  \int_{[0,T]\times\Omega}\textless \rho \sigma,\nabla_x \phi(x)\textgreater u(t)dtdx= 
  \int_{[0,T]\times\Omega}\rho\textless |\nabla_x G^{'}(\rho)|^{q(x)-2}\nabla_x G^{'}(\rho),\nabla_x\phi(x)\textgreater
  u(t) dtdx
.\end{equation}

Finally, we deduce that the sequence 
$div_x(\rho^h|\nabla_x G^{'}(\rho^h)|^{q(x)-2}\nabla_x G^{'}(\rho^h))$
converges to \\ $div_x(\rho |\nabla_x G^{'}(\rho)|^{q(x)-2}\nabla_x G^{'}(\rho))$ in $[C^{\infty}_c([0,T]\times\Omega)]^{'}$.
$\blacksquare$
\end{prf}
\end{prf}
\subsection{Existence and uniqueness of solution}
In this section, we show existence and uniqueness of weak solutions of the parabolic $q(x)$-Laplacian equation
\eqref{equation q(x)}.
\begin{thm}
Assume that $G$ satisfies \rm{$(A_1)$} and the initial datum $\rho_0$ satisfy \rm{$(A_3)$}.\\
Then, the sequence $(\rho^h)_h$  converges strongly to a positive function $\rho(t,x)$ and
 $\rho\in L^{\infty}([0,\infty[\times\Omega).$ Also $\rho$ is a weak solution of the equation
 \eqref{equation q(x)}, ie, for all $\phi(t,x)\in C^{\infty}_c([0,\infty[\times\Omega)$, 
 $supp\phi(.,x)\subset [-T,T]$, for $0<T<\infty,$ we have:
 \begin{equation}
\int_{[0,\infty[\times\Omega}\rho \left[\frac{\partial \phi(t,x)}{\partial t}-
\textless|\nabla_x G^{'}(\rho)|^{q(x)-2}\nabla_x G^{'}(\rho),\nabla_x\phi(t,x)\textgreater\right]dtdx=
-\int_{\Omega}\rho_0 \phi(0,x)dx 
. \end{equation}
\end{thm}
\begin{prf}
 Using \eqref{euler}:
 \begin{equation}\label{existence 1}
  \int_{[0,\infty[\times\Omega}\frac{\rho^k-\rho^{k-1}}{h}\phi(t,x)dtdx+
  \int_{[0,\infty[\times\Omega}\rho^h\textless |\nabla_x G^{'}(\rho^h)|^{q(x)-2}\nabla_x G^{'}(\rho^h),
  \nabla_x\phi (t,x)\textgreater dtdx=0(h)
, \end{equation}
where $0(h)$ tends to $0$ when $h$ tends to $0.$\\
Note that:
\begin{equation}
 \int_{[0,\infty[\times\Omega}\frac{\rho^k-\rho^{k-1}}{h}\phi(t,x)dtdx=
 \int_{[0,\infty[\times\Omega}\rho^h\left[ \frac{\phi(t-h,x)-\phi(t,x)}{h}\right]dtdx-
 \frac{1}{h}\int_0^h\int_{\Omega}\rho^h\phi(t-h,x)dtdx
.\end{equation}
 Replacing the previous relation in \eqref{existence 1}, we have:
\begin{eqnarray}\label{theo2}
 && \int_{[0,\infty[\times\Omega}\rho^h\left[ \frac{\phi(t-h,x)-\phi(t,x)}{h}\right]dtdx-
 \frac{1}{h}\int_0^h\int_{\Omega}\rho^h\phi(t-h)dtdx+\\\nonumber
&& \int_{[0,\infty[\times\Omega}\rho^h\textless |\nabla_xG^{'}(\rho^h)|^{q(x)-2}\nabla_x
G^{'}(\rho^h),\nabla_x\phi (t,x)\textgreater dtdx=0(h)
.\end{eqnarray}
 We tend $h$ to $0$ in \eqref{theo2} and  use theorem \eqref{theo1} to obtain:
 \begin{equation}
  \int_{[0,\infty[\times\Omega}\rho \left[\frac{\partial \phi(t,x)}{\partial t}-
\textless|\nabla_x G^{'}(\rho)|^{q(x)-2}\nabla_x G^{'}(\rho),\nabla_x\phi(t,x)\textgreater\right]dtdx=
-\int_{\Omega}\rho_0 \phi(0,x)dx 
 .\end{equation}
We conclude that $\rho$ is a weak solution of the parabolic $q(x)$-Laplacian equation \eqref{equation q(x)}.
$\blacksquare$
\end{prf}
\begin{thm}
 Assume that $p(.):\Omega\rightarrow ]1,+\infty[$ satisfy \rm{$(A_1)$}, $G$  satisfy \rm{$(A_2)$} and \rm{$(A_4)$};
 and $\rho_0$ satisfy
 \rm{($A_3)$}. Let $\rho^1$ and $\rho^2$ be two weak solutions of \eqref{equation q(x)} satisfying 
 $\frac{\partial\rho^{i}}{\partial t} \in L^{1}(\Omega)$, 
 for $i=1,2$,  with initial datum $\rho^1(0,.)$ and
 $\rho^2(0,.)$ respectively satisfying \rm{$(A_3)$}.\\
If $\rho^1(0,.)=\rho^2(0,.)$, then:
\begin{equation}
 \int_{\Omega}[\rho^1(T,x)-\rho^2(T,x)]^{+}dx\leq 0, 
\end{equation}
for all $T\geq0.$
\end{thm}
\begin{prf}
Define $\theta_{\delta}:\mathbb{R}\rightarrow [0,1]$, by:
\begin{eqnarray}\label{parabolic p,c}
\theta_{\delta}(s)=\left\{\begin{array}{rl}
 & 0\quad \mbox{if}\quad s\leq 0\\
 &\frac{s}{\delta} \quad \mbox{if}\quad 0\leq s\leq \delta\\
 & 1\quad \mbox{if} \quad s\geq \delta.
\end{array}\right.\quad
,\end{eqnarray}
 By using definition of the weak solution, we have:
 \begin{equation}\label{uniciterA}
  \int_{[0,T]\times\Omega}\phi\frac{\partial }{\partial t}(\rho^1(t,x)-\rho^2(t,x))dtdx=
  -\int_{[0,T]\times\Omega}\textless\rho^1|\nabla G^{'}(\rho^1)|^{q(x)-2}\nabla G^{'}(\rho^1)-
  \rho^2|\nabla G^{'}(\rho^2)|^{q(x)-2}\nabla G^{'}(\rho^2),\nabla \phi \textgreater dtdx.
 \end{equation}
We use $\theta_{\delta}(G^{'}(\rho^1)-G^{'}(\rho^2))$ in \eqref{uniciterA}; we have:
\begin{eqnarray}
 \int_{[0,T]\times\Omega}\theta_{\delta}(G^{'}(\rho^1)-G^{'}(\rho^2))\frac{\partial }{\partial t}(\rho^1(t,x)-
 \rho^2(t,x))dtdx
 = \\\nonumber
 -\int_{[0,T]\times\Omega}\textless\rho^1|\nabla G^{'}(\rho^1)|^{q(x)-2}\nabla G^{'}(\rho^1)-
 \rho^2|\nabla G^{'}(\rho^2)|^{q(x)-2}\nabla G^{'}(\rho^2),\nabla(\theta_{\delta}(G^{'}(\rho^1)-G^{'}(\rho^2)))
  \textgreater dtdx=\\
  -\frac{1}{\delta}\int_{\Omega_{T,\delta}}\rho^1\textless \nabla G^{'}(\rho^1)|^{q(x)-2}\nabla G^{'}(\rho^1)-
  |\nabla G^{'}(\rho^2)|^{q(x)-2}\nabla G^{'}(\rho^2),\nabla G^{'}(\rho^1)-\nabla G^{'}(\rho^2)\textgreater+\\\nonumber
  -\frac{1}{\delta}\int_{\Omega_{T,\delta}}\textless (\rho^1-\rho^2)|\nabla G^{'}(\rho^2)|^{q(x)-2}\nabla G^{'}(\rho^2),
  \nabla G^{'}(\rho^1)-\nabla G^{'}(\rho^2) \textgreater
\end{eqnarray}
Where $$ \Omega_{T,\delta}:=\Omega_T \cap \{0<G^{'}(\rho^1)-G^{'}(\rho^2)\leq \delta\}$$
and $\Omega_T:=[0,T]\times\Omega.$\\
Since $y\mapsto |y|^{q(x)-2}y$ is monotone, we have $$-\frac{1}{\delta}\int_{\Omega_{T,\delta}}\rho^1\textless
\nabla G^{'}(\rho^1)|^{q(x)-2}\nabla G^{'}(\rho^1)-
  |\nabla G^{'}(\rho^2)|^{q(x)-2}\nabla G^{'}(\rho^2),\nabla G^{'}(\rho^1)-\nabla G^{'}(\rho^2)\textgreater\leq 0.$$ 
Furthermore, on $\Omega_{T,\delta}$:
$$|\rho^1-\rho^2|=|(G^*)^{'}(G^{'}(\rho^1))-(G^*)^{'}(G^{'}(\rho^2))|
\leq \delta \sup_{s\in [0,G^{'}(M_2)]}(G^*)^{''}(s)$$  
Then,
$$-\frac{1}{\delta}\int_{\Omega_{T,\delta}}\textless (\rho^1-\rho^2)|\nabla G^{'}(\rho^2)|^{q(x)-2}\nabla G^{'}(\rho^2),
  \nabla G^{'}(\rho^1)-\nabla G^{'}(\rho^2) \textgreater  \precsim |\Omega_{T,\delta}|.$$
  If $\delta \rightarrow 0^{+}$, then $|\Omega_{T,\delta}|\rightarrow 0$ and 
  $\theta_{\delta}(G^{'}(\rho^1)-G^{'}(\rho^2))\rightarrow sign^{+}(G^{'}(\rho^1)-G^{'}(\rho^2))=sign^{+}(\rho^1-\rho^2);$
  with
  $sign(s)=\frac{s}{|s|}$ for all $s\in\mathbb{R}^*.$
  Then, 
  \begin{equation}
\int_{[0,T]\times\Omega}\frac{\partial(\rho^1-\rho^2)^{+}}{\partial t}=
\int_{[0,T]\times\Omega}sign^{+}(\rho^1-\rho^2)\frac{\partial }{\partial t}(\rho^1-\rho^2)\leq 0.
  \end{equation}
  This implies
  \begin{equation}
   \int_{\Omega}(\rho^1(T,x)-\rho^2(T,x))^{+}dx\leq 0 
  \end{equation}
for all $T\geq 0.$
Then the solution of $q(x)$- Laplacian equation  \eqref{equation q(x)}  is unique.
$\blacksquare$
\end{prf}
  
Aboubacar Marcos and Ambroise Soglo\\
Institut de Math\'ematiques et de Sciences Physiques\\
Universit\'e d'Abomey-Calavi\\
E-mail: abmarcos24@gmail.com or abmarcos@imsp-uac.org\\
E-mail: ambroise.soglo@yahoo.fr

\begin{thebibliography}{9}
 \bibitem{aguey}M.Agueh, \emph{Finsler structure in $p$- Wasserstein space and gradient flows}, 
 university of Victoria, P.O.Box
  3060, STNCSC, Victoria, BC, VBW3R4, Canada.
  \bibitem{agueh}M. Agueh, \emph{Existence of solutions to degenerate parabolic equations via
the Monge-Kantorovich theory}, Adv. Differential Equations, 10 (2005), 309-360.
\bibitem{ambrosio} L. Ambrosio, N. Gigli, G. Savaré,\emph{Gradient flow in Metric Spaces and in the space of probability
Measures}, Lectures in Mathematics, Birkhauser,2005.
\bibitem{brezis}H. Brezis, \emph{Analyse fonctionnelle théorie et applications}, 
Masson Paris New York Barcelone Milan Mexico Sao Paulo 1987.
\bibitem{chen} Y.M. Chen.S.Levine, M.Ra,\emph{Variable exponent, linear growth functionals in image restoration},SIAM J.Appl.
Math; 66 (2006); 1383-1406.
 \bibitem{mawhin} G. Dinca and J. Mawhin,\emph{variational and topological methods for Dirichlet problems with p-Laplacian}, 
  Portugaliae Mathematica vol.58 Fasc.3-2001.
  \bibitem{dinca2} G. Dinca, P. Matei \emph{Geometry of Sobolev spaces with variable exponent: Smoothness and 
  uniform convexity} C.R. Acad. Sci. Paris, Ser. 1347 (2009) 885-889.
  \bibitem{fan zao} Fan and Zhao, \emph{On the spaces $L^{p(x)}(\Omega)$ and $W^{m,p(x)}(\Omega)$}, Department of mathematics,
  Lanzhou university, Lanzhou 730000.
\bibitem{figalli} A. Figalli, W. Gangbo and T. Yolcu, \emph{ Parabolic equations involving a Lagrangian}
. Ann. Scuola Norm. Sup. Pisa Cl. Sci. (5), 10, no. 1, 207-252, 2011.
\bibitem{ekeland} D.G. De Figueiredo,\emph{Lectures on The Ekeland Variational Principle with Applications and Detours},
Tata Institute of Fundamental Research, Bombay 1989
\bibitem{gangbo} W. Gangbo, R. McCann\emph{ The geometry of optimal  transport}.  Acta Mathematica, Vol 177, no. 2, 113--161, 
1996.
 \bibitem{gangboo} W. Gangbo, R. McCann \emph{ Optimal Maps in Monge's mass transport problem}.
C.R. Acad. Sci. Paris, t.321, S\'erie I, 1653-1658, 1995. 
   \bibitem{jordan}R. Jordan, D. Kinderlehrer and F. Otto, \emph{The variational formulation of the Fokker-Planck
equation}, SIAM J. Math. Anal., 29 (1998), 1-17.
\bibitem{Huashui} Huashui Zhan,\emph{The stability of evolutionary p(x)-Laplacian equation}, Zhan Boundary value Problems 
(2017).
\bibitem{Lian} Lian,S,Gao,W, Yan,H,H,Cao,C:\emph{Existence of solutions to an Dirichlet problem of evolutional 
p(x)-Laplacian equations},Ann.inst.Henri Poincar\'e, Anal. Non lin\'eaire 29,377-399 (2012).
\bibitem{Livre variable} Lars Diening,Petteri Harjulehto, Peter Hasto, Michael R\r{u}\v{z}i\u{c}ka, 
\emph{Lebesgue and Sobolev spaces with variable
exponents}, Lecture Notes in Mathematics 2017.
\bibitem{Myers} T.G Myers,\emph{Thin filmswith light surface tension}, SiAM Review, 40 (1998); 441-462.
\bibitem{rund} H.Rund,\emph{The Differential Geometry of Finsler Spaces, Springer},Verlag 1959.
\bibitem{Ruzicka} M. R\^{u}\v{z}icka; \emph{Electrorheological fluids: modeling and Mathematical theory},Lecture notes
in Mathematics,
1748, Springer-Verlag,Berlin,2000.
\bibitem{Yamna}G.Yamna and H.El Ouardi,\emph{Study of solutions to a class of certain parabolic systems with variable 
exponents},
EDP sciences, 2017.
\bibitem{Wen} H. Zhan and Jie Wen,
\emph{Evolutionary p(x)-Laplacian Equation Free From the limitation of the bounded value}, 2016, Texas state
university.
\bibitem{Huashui}  Huashui Zhan \emph{The stability of evolutionary $p(x)$-Laplacian equation}, Zhan
Boundary Value Problems
 (2017) 
\bibitem{Zhikov} V.V. Zhikov, S.M. Kozlov, O.A. Oleinik,
\emph{Homogenization of differential operators and integral functionals},
Translated from the Russian by G.A. Yosifian.Springer-Verlag, Berlin, 1994.
\bibitem{Zhan H} Zhan.H: \emph{The boundary value condition of an evolutionary p(x)-Laplacian equation}.
Bound value Probl.2015, 112 (2015).
 \end{thebibliography}
\end{document}